\documentclass[final]{article}

\usepackage{graphicx}
\usepackage{authblk}
\usepackage{amssymb}
\usepackage{epstopdf}
\usepackage[sans]{dsfont}
\usepackage[T1]{fontenc}
\usepackage[english]{babel}
\usepackage{latexsym}
\usepackage{subfigure} 
\usepackage{mathrsfs}
\usepackage{amscd}
\usepackage{color}
\usepackage{float}
\usepackage{bm}
\usepackage{amsmath}
\usepackage{amsfonts}
\usepackage{tikz}
\usetikzlibrary{calc,trees,positioning,arrows,chains,shapes.geometric,%
    decorations.pathreplacing,decorations.pathmorphing,shapes,%
    matrix,shapes.symbols}
\tikzset{
>=stealth',
  punktchain/.style={
    rectangle, 
    rounded corners, 
     fill=cyan!40,
    draw=black, very thick,
    text width=10em, 
    minimum height=2em, 
    text centered, 
    on chain},
  line/.style={draw, thick, <-},
  element/.style={
    tape,
    top color=white,
    bottom color=blue!50!black!60!,
    minimum width=8em,
    draw=blue!40!black!90, very thick,
    text width=10em, 
    minimum height=2.5em, 
    text centered, 
    on chain},
  every join/.style={->, thick,shorten >=1pt},
  decoration={brace},
  tuborg/.style={decorate},
  tubnode/.style={midway, right=2pt},
}
\numberwithin{equation}{section}
	\usepackage{enumerate}
\usepackage{amsthm}
\usepackage[numbers]{natbib}
\usepackage[bookmarks=true,colorlinks=true,linkcolor={blue},urlcolor={blue}, citecolor={blue},pdfstartview={XYZ null null 1.22}]{hyperref}%

\def\vm{{\mathcal{V}}}
\newtheorem{remark}{Remark}

\def\RR{\mathbb R}

\def\be{\begin{equation}}
\def\ee{\end{equation}}
\def\bea{\begin{eqnarray}}
\def\eea{\end{eqnarray}}

\begin{document}
\title{Uncertainty Quantification in Control Problems for Flocking Models}

\author{Giacomo Albi\thanks{Fakult\"at f\"ur Mathematik, Technische Univarsit\"at M\"unchen, Germany }\qquad  Lorenzo Pareschi\thanks{Department of Mathematics and Computer Science, University of
Ferrara, Italy}
\qquad 
Mattia Zanella\thanks{Department of Mathematics and Computer Science, University of
Ferrara, Italy}}

\maketitle

\begin{abstract}
In this paper the optimal control of flocking models with random inputs is investigated from a numerical point of view. The effect of uncertainty in the interaction parameters is studied for a Cucker-Smale type model using a generalized polynomial chaos (gPC) approach. Numerical evidence of threshold effects in the alignment dynamic due to the random parameters is given. The use of a selective model predictive control permits to steer the system towards the desired state even in unstable regimes. 
\end{abstract}


 
\tableofcontents

\section{Introduction}

The aggregate motion of a multi-agent system is frequently seen in the real world. Common examples are represented by schools of fishes, swarms of bees and herds of sheep, each of them natural phenomena with important applications in many fields such as biology, engineering and economy \cite{MR2761862}. As a consequence, the significance of new mathematical models, for understanding and predicting these complex dynamics, is widely recognized. Several heuristic rules of flocking have been introduced as alignment, separation and cohesion \cite{R,V}. Nowadays these mathematical problems, and their constrained versions, are deeply studied both from the microscopic viewpoint \cite{ABC,MR2596552,ZZZ} as well as their kinetic and mean-field approximations \cite{AHP, APZ,MR2324245,DLR14,DHL}. We refer to \cite{MR2761862} for a recent introduction on the subject. 

An essential step in the study of such models is represented by the introduction of stochastic parameters reflecting the uncertainty in the terms defining the interaction rules. In \cite{AH10,HLL09,YEE2009} the authors are concerned with the study of self-organized system including noise term modeling the interaction with the surrounding environment. In this paper we focus on the case where the uncertainty acts directly in the parameter characterizing the interaction dynamic between the agents.

We present a numerical approach having roots in the numerical techniques for Uncertainty Quantification (UQ) and Model Predictive Control (MPC). Among the most popular methods for UQ, the generalized polynomial chaos (gPC) has recently received deepest attentions \cite{X}. Jointly with Stochastic Galerkin (SG) this class of numerical methods are usually applied in physical and engineering problems, for which fast convergence is needed. Applications of gPC-Galerkin schemes to flocking dynamics, and their controlled versions, is almost unexplored in the actual state of art. 

We give numerical evidence of threshold effects in the alignment dynamic due to the random parameters. In particular the presence of a negative tail in the distribution of the random inputs lead to the divergence of the expected values for the system velocities.
The use of a selective model predictive control permits to steer the system towards the desired state even in such unstable regimes. 

The rest of the article is organized as follows. In Section \ref{sec:CS} we introduce briefly a Cucker-Smale dynamic with interaction function depending on stochastic parameters and analyze the system behavior in the case of uniform interactions. The gPC approach is then summarized in Section \ref{sec:gPC}. Subsequently, in Section \ref{sec:UQ} we consider the gPC scheme in a constrained setting and derive a selective model predictive approximation of the system. Next, in Section \ref{sec:Numeric} we report several numerical experiments which illustrate the different features of the numerical method. Extensions of the present approach are finally discussed in Section \ref{sec:Conc}.

\section{Cucker-Smale dynamic with random inputs}\label{sec:CS}
	We introduce a Cucker-Smale type \cite{MR2324245} differential system depending on a random variable $\theta\in \Omega\subseteq \mathbb{R}$ with a given distribution $f(\theta)$. Let $(x_i,v_i)\in\mathbb{R}^{2d}, d\geq 1$, evolving as follows
\begin{equation}\begin{cases}
\dot{x}_i(\theta,t)=v_i(\theta,t) \\
\dot{v}_i(\theta,t)=\dfrac{K(\theta,t)}{N}\displaystyle\sum_{j=1}^N H(x_i,x_j)(v_j(\theta,t)-v_i(\theta,t))
\end{cases}
\label{eq:CS}
\end{equation}
where $K$ is a time dependent random function characterizing the uncertainty in the interaction rates and $H(\cdot,\cdot)$ is a symmetric function describing the dependence of the alignment dynamic from the agents positions. A classical choice of space dependent interaction function is related to the distance between two agents 
\begin{equation}
H(x,y)=\dfrac{1}{(1+|x-y|^2)^\gamma},
\label{eq:interaction_space}
\end{equation} 
where $\gamma\ge0$ is a given parameter. {Mathematical results concerning the system behavior in the deterministic case ($K\equiv 1$) can be found in \cite{MR2324245}. In particular unconditional alignment emerges for $\gamma<1/2$. Let us observe that, even for the model with random inputs \eqref{eq:CS}, the mean velocity of the system is conserved in time
\begin{equation}
\vm (\theta,t)=\frac1{N}\sum_{i=1}^N {v}_i(\theta,t),\qquad \frac{d}{dt} \vm (\theta,t) = 0,
\end{equation}
since the symmetry of $H$ implies
\[
\sum_{i,j=1}^N H(x_i,x_j)v_j(\theta,t)=\sum_{i,j=1}^N H(x_i,x_j)v_i(\theta,t).
\]
Therefore we have $\vm (\theta,t)=\vm (\theta,0)$.
}

\subsection{The uniform interaction case}
\label{uniform}
To better understand the leading dynamic let us consider the simpler uniform interaction case when $H\equiv 1$, leading to the following equation for the velocities 
\begin{equation}
\dot{v}_i(\theta,t)=\dfrac{K(\theta,t)}{N}\displaystyle\sum_{j=1}^N(v_j(\theta,t)-v_i(\theta,t))=K(\theta,t)(\vm (\theta,0)-v_i(\theta,t)).
\label{eq:CShomo}
\end{equation}
The differential equation \eqref{eq:CShomo} admits an exact solution depending on the random input $\theta$. More precisely if the initial velocities are deterministically known we have that 
\begin{equation}\label{eq:exact_homo}
v_i(\theta,t)=\vm +(v_i(0)-\vm ) \exp{\Big\{-\int_0^t K(\theta,s)ds\Big\}},
\end{equation}
where $\vm =\vm (0)$ is the mean velocity  of the system.
In what follows we analyze the evolution of \eqref{eq:exact_homo} for different choices of $K(\theta,t)$ and of the distribution of the random variable $\theta$. 
\subsubsection*{Example 1}
Let us consider a random scattering rate written in terms of the following decomposition 
\begin{equation}\label{eq:K_normal}
K(\theta,t)=k(\theta)h(t)
\end{equation}
where $h(t)$ is a nonnegative function depending on $t\in\RR^+$. The expected velocity of the $i$-th agent is defined by 
\begin{equation}\label{eq:expected_velocity}
\bar{v}_i(t)=\mathbb{E}_{\theta}[v_i(\theta,t)]=\int_{\Omega}v_i(\theta,t)f(\theta)d\theta
\end{equation} 
then each agent evolves its expected velocity according to
\begin{equation}\label{eq:expectation}
\bar{v}_i(t)=\int_{\Omega}\left[ \vm +(v_i(0)-\vm ) \exp{\Big\{-k(\theta)\int_0^t h(s)ds\Big\}}\right]f(\theta)d\theta.
\end{equation}
For example, let us chose $k(\theta)=\theta$, where the random variable is normally distributed, i.e. $\theta\sim\mathcal{N}(\mu,\sigma^2)$. Then, for each $i=1,\dots,N$, we need to evaluate the following integral 
\begin{equation}
\vm +\dfrac{v_i(0)-\vm }{\sqrt{2\pi\sigma^2}}\int_{\RR}\exp{\Big\{-\theta\int_0^t h(s)ds\Big\}}\exp{\Big\{-\dfrac{(\theta-\mu)^2}{2\sigma^2}\Big\}}d\theta.
\end{equation}
The explicit form is easily found through standard techniques and yields
\begin{equation}\label{homo_normal}
\bar{v}_i(t)=\vm +(v_i(0)-\vm )\exp{\Big\{ -\mu\int_0^t h(s)ds+\dfrac{\sigma^2}{2}\left(\int_0^t h(s)ds\right)^2 \Big\}}.
\end{equation}
From \eqref{homo_normal} we observe a threshold effect in the asymptotic convergence of the mean velocity of each agent toward $\vm $. It is immediately seen that if
\begin{equation}
\int_0^t h(s)ds>\dfrac{2\mu}{\sigma^2}
\end{equation}
it follows that, for $t\rightarrow+\infty$, the expected velocity $\bar{v}_i$ diverges. {In particular, if $h(s)\equiv 1$ we have that the solution starts to diverge as soon as $t>\mu/\sigma^2$. Note that, this threshold effect is essentially due to the negative tail of the normal distribution. In fact, if we now consider a random variable taking only nonnegative values, for example exponentially distributed $\theta\sim\mathrm{Exp}(\lambda)$ for some positive parameter $\lambda>0$, from equation \eqref{eq:expectation} we obtain 
\begin{equation}
\bar{v}_i(t)=\vm +(v_i(0)-\vm )\int_0^{+\infty} e^{-\theta t} \lambda e^{-\lambda\theta}d\theta,
\end{equation}
which corresponds to 
\begin{equation}
\bar{v}_i(t)=\vm +(v_i(0)-\vm )\dfrac{\lambda}{t+\lambda},
\end{equation}
and therefore $\bar{v}_i(t)\to \vm $ as $t\to \infty$. Then independently from the choice of the rate $\lambda>0$ we obtain for each agent convergences toward the average initial velocity of the system. Finally, in case of a uniform random variable $\theta\sim U([a,b])$ we obtain 
\begin{equation}
\bar{v}_i(t)=\vm +(v_i(0)-\vm )\int_a^b \dfrac{1}{b-a} e^{-\theta t}d\theta
\end{equation}
that is 
\begin{equation}
\bar{v}_i(t)=\vm +\dfrac{v_i(0)-\vm }{b-a}\left(\dfrac{e^{-a t}}{t}-\dfrac{e^{-bt}}{t}\right),
\end{equation}
which implies the divergence of the system in time as soon as $a$ assumes negative values.
 
\subsubsection*{Example 2}\label{rem:variance}
Next we consider a random scattering rate with time-dependent distribution function, that is 
\begin{equation}
K(\theta,t)=\theta(t)
\end{equation} 
with $\theta(t)\sim f(\theta,t)$. As an example we investigate the case of a normally distributed random parameter with given mean and time dependent variance,  $\theta\sim\mathcal{N}(\mu,\sigma^2(t))$. It is straightforward to rewrite $\theta$ as a translation of a standard normal-distributed variable $\tilde{\theta}$, that is
\begin{align}
\theta=\mu+\sigma(t)\tilde{\theta},
\end{align}
where $\tilde{\theta}\sim\mathcal{N}(0,1)$. The expected velocities read 
\begin{equation}
\bar{v}_i(t)=\vm +\frac{(v_i(0)-\vm )}{\sqrt{2\pi}}\int_{\RR} \exp{\Big\{-\mu t-\tilde{\theta}\int_0^t \sigma(s)ds\Big\}}\exp\left\{-\frac{\tilde{\theta}^2}{2}\right\} d\tilde{\theta},
\end{equation}
which correspond to 
\begin{equation}
\bar{v}_i(t)=\vm +(v_i(0)-\vm )\exp{\Big\{-\mu t +\dfrac{1}{2}\left(\int_0^t\sigma(s)ds\right)^2\Big\}}.
\end{equation}
Similarly to the case of a time independent normal variable a threshold effect occurs for large times, i.e. the following condition on the variance of the distribution 
\begin{equation}\label{eq:resh2}
\left(\int_0^t\sigma(s)ds\right)^2>{2\mu t}
\end{equation}
implies the divergence of the system \eqref{eq:CShomo}.
As a consequence instability can be avoided by assuming a variance decreasing sufficiently fast in time. The simplest choice is represented by $\sigma(t)=1/t^{\alpha}$ for some $\alpha\in[1/2,1)$. The condition \eqref{eq:resh2} becomes
\begin{equation}
\left(\dfrac{t^{1-\alpha}}{1-\alpha}\right)^2>{2\mu t}.
\end{equation}
For example, if $\alpha=1/2$ the previous condition implies that the system diverges for each $\mu<2$.

\section{A gPC based numerical approach}
\label{sec:gPC}
In this section we approximate the Cucker-Smale model with random inputs using a generalized polynomial chaos approach. For the sake of clarity we first recall some basic facts concerning gPC approximations.
\subsection{Preliminaries on gPC approximations}
Let $(\Omega,\mathcal{F},P)$ be a probability space, that is a ordered triple with $\Omega$ any set, $\mathcal{F}$ a $\sigma-$algebra and $P:\mathcal{F}\rightarrow [0,1]$ a probability measure on $\mathcal{F}$, where we define a random variable
\[
\theta:(\Omega,\mathcal{F})\rightarrow(\mathbb{R},\mathcal{B}_{\mathbb{R}}),
\]
with $\mathcal{B}_{\mathbb{R}}$ the Borel set of $\mathbb{R}$. Moreover let us consider $S\subset \mathbb{R}^d,d\ge 1$ and $[0,T]\subset\mathbb{R}$ certain spatial and temporal subsets. For the sake of simplicity we focus on real valued functions depending on a single random input
\begin{equation}\label{eq:initial_gPC}
g(x,\theta,t): S\times \Omega\times T\rightarrow \mathbb{R}^d, \qquad g\in L^2(\Omega,\mathcal{F},P).
\end{equation}
In any case it is possible to extend the set-up of the problem to a $p-$dimensional vector of random variables $\bm{\theta}=(\theta_1,\dots,\theta_p)$, see \cite{Vos}.  Let us consider now the linear space of polynomials of $\theta$ of degree up to $M$, namely $\mathbb{P}_M(\theta)$. From classical results in approximation theory it is possible to represent the distribution of random functions with orthogonal polynomials $\{\Phi_k(\theta)\}_{k=0}^M$, i.e. an orthogonal basis of $L^2(\Omega,\mathcal{F},P)$
\[
\mathbb{E}_\theta[ \Phi_h(\theta) \Phi_k(\theta)]=\mathbb{E}_{\theta}[\Phi_h(\theta)]^2 \delta_{hk}
\]
with $\delta_{hk}$ the Kronecker delta function. Assuming that the probability law $P(\theta^{-1}(B)), \forall B\in\mathcal{B}_{\mathbb{R}}$, involved in the definition of the introduced function  $g(x,\theta,t)$ has finite second order moment, then the complete polynomial chaos expansion of $g$ is given by 
\begin{equation}
g(x,\theta,t)=\sum_{m\in\mathbb{N}} \hat{g}_{m}(x,t)\Phi_{m}(\theta).
\end{equation}
Accordingly to the Askey-scheme, result that pave a connection between random variables and orthogonal polynomials \cite{X,XK}, we chose a set of polynomials which constitutes the optimal basis with respect to the distribution of the introduced random variable in agreement with Table \ref{tab:pol}.
\begin{table}[t]
\caption{The different gPC choices for the polynomial expansions}
\begin{center}
\begin{tabular}{c | c |  c}\label{tab:Ask}
Probability law of $\theta$ & Expansion polynomials & Support \\ \hline\hline
Gaussian & Hermite    & $(-\infty,+\infty)$ \\
Uniform    & Legendre & $[a,b]$\\
Beta         & Jacobi     & $[a,b]$\\
Gamma    & Laguerre & $[0,+\infty)$\\
Poisson    & Charlier   & $\mathbb{N}$\\
\end{tabular}
\end{center}
\label{tab:pol}
\end{table}

Let us consider now a general formulation for a randomly perturbed problem 
\begin{equation}\label{eq:problem}
\mathcal{D}(x,t,\theta;g)=f(x,t,\theta)
\end{equation}
where we indicated with $\mathcal{D}$ a differential operator. In general the randomness introduced in the problem by $\theta$ acts as a perturbation of $\mathcal{D}$, of the function $g$ or occurs as uncertainty of the initial conditions. In this work we focus on the first two aspects assuming that initial positions and velocities are deterministically known.

The generalized polynomial chaos method approximate the solution $g(x,\theta,t)$ of \eqref{eq:problem} with its $M$th order polynomial chaos expansion and considers the Galerkin projections of the introduced differential problem, that is
\begin{equation}\label{eq:problem_projection}
\mathbb{E}_{\theta}\left[\mathcal{D}(x,t,\theta;g) \cdot \Phi_h(\theta)\right]=\mathbb{E}_{\theta}\left[f(x,\theta,t)\cdot \Phi_h(\theta)\right],\qquad h=0,1,\dots,M.
\end{equation}
Due to the Galerkin orthogonality between the linear space $\mathbb{P}_M$ and the error produced in the representation of $g(x,\theta,t)$ with a truncated series, it follows that from \eqref{eq:problem_projection} we obtain a set of $M+1$ purely deterministic equations for the expansion coefficients $\hat{g}_m(x,t)$. These subproblems can be solved through classical discretization techniques. From the numerical point of view through a gPC-type method it is possible to achieve an exponential order of convergence to the exact solution of the problem, unlike Monte Carlo techniques for which the order is $O(1/\sqrt{M})$ where $M$ is the number of samples.


\subsection{Approximation gPC of the alignment model}
We apply the described gPC decomposition to the solution of the non-homogeneous differential equation $v_i(\theta,t)$ in \eqref{eq:exact_homo} and to the stochastic scattering rate $K(\theta,t)$, i.e.  
\begin{equation}
v_i^M(\theta,t)=\sum_{m=0}^M \hat{v}_{i,m}(t)\Phi_m(\theta), \qquad K^M(\theta,t)=\sum_{l=0}^M \hat{K}_l(t)\Phi_l(\theta)
\end{equation}
where
\begin{align}
\hat{v}_{i,m}(t)=\mathbb{E}_{\theta}\left[v_i(\theta,t)\Phi_m(\theta)\right] \qquad \hat{K}_l(t)=\mathbb{E}_{\theta}\left[K(\theta,t)\Phi_l(\theta)\right],
\end{align}
 we obtain the following polynomial chaos expansion
 \begin{equation*}
 \begin{aligned}
\dfrac{d}{dt}\sum_{m=0}^M \hat{v}_{i,m}\Phi_m(\theta)=&\dfrac{1}{N}\sum_{j=1}^N H(x_i,x_j) \sum_{l,m=0}^M \hat{K}_l(t)(\hat{v}_{j,m}-\hat{v}_{i,m})\Phi_l(\theta)\Phi_m(\theta).\\
\end{aligned}
\end{equation*} 
Multiplying the above expression by an orthogonal element of the basis $\Phi_h(\theta)$ and integrating with respect to the distribution of $\theta$ 
\begin{equation*}
\begin{aligned}
\mathbb{E}_{\theta}&\left[\sum_{m=0}^M\dfrac{d}{dt}\hat{v}_{i,m}\Phi_m(\theta)\Phi_h(\theta)\right] \\
 &\qquad =\mathbb{E}_{\theta}\left[\dfrac{1}{N}\sum_{j=1}^N H(x_i,x_j) \sum_{l,m=0}^M \hat{K}_l (t)(\hat{v}_{j,m}-\hat{v}_{i,m})\Phi_l(\theta)\Phi_m(\theta)\Phi_h(\theta)\right]\\
\end{aligned}
\end{equation*}
we find an explicit system of ODEs 
\begin{equation}\label{eq:decomposition} 
\begin{aligned}
\dfrac{d}{dt}\hat{v}_{i,h}(t)=&\dfrac{1}{\|\Phi_h \|^2 N}\sum_{j=1}^N H(x_i,x_j)\sum_{m=0}^M (\hat{v}_{j,m}-\hat{v}_{i,m})\sum_{l=0}^M \hat{K}_l(t) e_{lmh}\\
=& \dfrac{1}{N}\sum_{j=1}^N H(x_i,x_j)\sum_{m=0}^M (\hat{v}_{j,m}-\hat{v}_{i,m})\hat{K}_{mh}(t),
\end{aligned}
\end{equation}
where $e_{lmh}=\mathbb{E}_{\theta}[\Phi_l(\theta)\Phi_m(\theta)\Phi_h(\theta)]$ and
$$\hat{K}_{mh}(t)= \dfrac{1}{\|\Phi_h \|^2}\sum_{l=0}^M \hat{K}_l(t) e_{lmh}.$$
We recall that the gPC numerical approach preserves the mean velocity of the alignment model \eqref{eq:CShomo}. In fact, from \eqref{eq:decomposition} follows
\begin{equation}\begin{split}
\sum_{i=1}^N \dfrac{d}{dt}\hat{v}_{i,h}(t)=&\dfrac{1}{N}\sum_{j,i=1}^N H(x_i,x_j)\sum_{m=0}^M \hat{K}_{mh}(t) \hat{v}_{j,m}\\
& -\dfrac{1}{N}\sum_{j,i=1}^N H(x_i,x_j)\sum_{m=0}^M \hat{K}_{mh}(t) \hat{v}_{i,m}=0,
\end{split}\end{equation}
thanks to the symmetry of $H$. More  generally it can be shown that if $$\frac{1}{N}\sum_{i=1}^Nv_i(\theta,t)=\vm ,$$ where $\vm $ is time-independent, then its gPC decomposition is also mean-preserving since
\begin{equation*}\begin{split}
\dfrac{1}{N}\sum_{i=1}^N\sum_{m=0}^M\mathbb{E}_{\theta}\left[v_i(\theta,t)\Phi_m(\theta)\right]\Phi_m(\theta)&=\sum_{m=0}^M \mathbb{E}_{\theta}\left[\dfrac{1}{N}\sum_{i=1}^N v_i(\theta,t)\Phi_m(\theta)\right]\Phi_m(\theta) \\
&=\vm \sum_{m=0}^M \mathbb{E}_{\theta}\left[1\cdot\Phi_m(\theta)\right]\Phi_m(\theta)
=\vm .
\end{split}\end{equation*}
\begin{remark}
The gPC approximation (\ref{eq:decomposition}) can be derived equivalently without expanding the kernel function $K(\theta,t)$. In this way one obtains 
\begin{equation}\label{eq:decomposition2} 
\begin{aligned}
\dfrac{d}{dt}\hat{v}_{i,h}(t)=&\dfrac{1}{N}\sum_{j=1}^N H(x_i,x_j)\sum_{m=0}^M (\hat{v}_{j,m}-\hat{v}_{i,m}){\widehat K}_{mh}\\
\end{aligned}
\end{equation} 
where now $${\widehat K}_{mh}(t)=\dfrac{1}{\|\Phi_h \|^2}\mathbb{E}_\theta[K(\theta,t)\Phi_m\Phi_h].$$
Note that, since in general $N \gg M$, the overall computational cost is $O(MN^2)$.
\end{remark}

\section{Selective control of the gPC approximation}\label{sec:UQ}
In order stabilize the gPC approximation of the Cucker-Smale type model \eqref{eq:CS} with random inputs, we introduce an additional term which acts as control of the approximated dynamic. 
More specifically we modify the approximation of the alignment model \eqref{eq:CS} by introducing a control term $\hat{u}_h$ to the $h$th component of its gPC approximation 
\begin{equation}\begin{split}
\dfrac{d}{dt}\hat{v}_{i,h}(t)=&\dfrac{1}{N}\displaystyle\sum_{j=1}^N H(x_i,x_j)\sum_{m=0}^M \hat{K}_{mh} (t)(\hat{v}_{j,m}(t)-\hat{v}_{i,m}(t))+\hat{u}_hQ(\hat{v}_{i,h})
\label{eq:CScontrol}
\end{split}\end{equation}
where $\hat{u}_h$ is a solution of 
\begin{equation}
\hat{u}_h=\textrm{arg}\min_{\hat{u}_h\in\RR}\left[\dfrac{1}{2}\int_0^T\dfrac{1}{N}\sum_{i=1}^N(\hat{v}_{i,h}(t)-\hat{v}_{d,h})^2dt + \dfrac{\nu}{2}\int_0^T\hat{u}_h(t)^2dt \right],
\label{eq:control}
\end{equation}
where $\nu>0$ is a regularization parameter and $(\hat{v}_{d,0},\hat{v}_{d,1},\dots,\hat{v}_{d,M})$ are the desired values for the gPC coefficients. For example 
\begin{equation}
\hat{v}_{d,h}=\mathbb{E}_{\theta}[v_d\Phi_h(\theta)]=v_d\mathbb{E}[\Phi_h(\theta)]=
\begin{aligned}\begin{cases}
 v_d & h=0 \\
 0     & h=1,\dots,M, 
\end{cases}\end{aligned}
\end{equation}
where $v_{d}$ is a desired velocity.

  Moreover the controller action is weighted by a function, $Q(\cdot)$, such that $Q(\hat{v}_{i,h})\in[-L,L],L>0$. Due to the dependence of the controller effect from the single agent velocity, we refer to this approach as selective control, see \cite{AP15}.  

 In order to tackle numerically the above problem, whose direct solution is prohibitively expansive for large numbers of individuals, we make use of model predictive control (MPC) techniques, also referred to as receding horizon strategy or instantaneous control \cite{MichalskaMayne1993aa}. These techniques has been used in \cite{AHP, AP15, APZ} in the case of deterministic alignment systems. 
 
\subsection{Selective model predictive control} 
The basic idea is to consider a piecewise constant control, 
$$\hat{u}_h(t)=\sum_{n=0}^{m-1}\hat{\bar{u}}_h^{n}\chi_{[t_n,t_{n+1}]}(t),$$ 
on a suitable time discretization.
In this way is possible to determine the  value of the control $\hat{u}_h^n \in \mathbb{R}$, solving for a state $\bar{\hat{v}}_{i,h}$ the (reduced) optimization problem                          
\begin{equation}\label{eq:pbic1}
\begin{aligned}
& \dfrac{d}{dt}\hat{v}_{i,h}(t)=\dfrac{1}{N}\displaystyle\sum_{j=1}^N H(x_i,x_j)\sum_{m=0}^M \hat{K}_{mh}(t) (\hat{v}_{j,m}(t)-\hat{v}_{i,m}(t))+\hat{u}_hQ(\hat{v}_{i,h}(t)) \\
& \hat{v}_{i,h}(t^n)=\bar{\hat{v}}_{i,h}, \\
& \hat{u}_h^n=\textrm{arg}\min_{\hat{u}_h\in\RR} \int_{t^n}^{t^{n+1}} \frac{1}{N}\sum_{i=1}^N\left(\frac{1}{2}(\hat{v}_{i,h}(t)-\hat{v}_{d,h})^2+ \frac{\nu}{2} \hat{u}_h(t)^2\right)ds, \hat{u}_h^n \in[\hat{u}_{h,L}^n,\hat{u}_{h,R}^n].
\end{aligned}
\end{equation}
Given the control  $\hat{u}_h^n$ on the time interval $[t^n, t^{n+1}]$, we let evolve $\hat{v}_{i,h}$ according to the dynamics
\begin{equation}\begin{split}\label{eq:controldynamics1}
\dfrac{d}{dt}\hat{v}_{i,h}=&\dfrac{1}{N}\displaystyle\sum_{j=1}^N H(x_i,x_j)\sum_{m=0}^M \hat{K}_{mh}(t) (\hat{v}_{j,m}(t)-\hat{v}_{i,m}(t))+\hat{u}_h^nQ(\hat{v}^n_{i,h}(t))
\end{split}\end{equation}
in order to obtain  the new state $\bar{\hat{v}}_{i,h} = \hat{v}_{i,h}(t^{n+1}).$ 
We again solve \eqref{eq:pbic1} to obtain $\hat{u}_h^{n+1}$ with the modified initial data and we repeat this procedure until we reach $n \Delta t = T.$   

The reduced optimization problem implies a reduction of the complexity of the initial problem since  to an optimization problem in a single real--valued variable $\hat{u}_h^n$. On the other hand the price to pay is that in general the solution of the problem is suboptimal respect to the full one  \eqref{eq:CScontrol}-\eqref{eq:control}.

The quadratic cost and a suitable discretization of \eqref{eq:controldynamics1} allows an explicit representation of $\hat{u}_h^n$ in terms of $\bar{\hat{v}}_{i,h}$ and $\hat{v}_{i,h}^{n+1}$, as a feedback controlled system as follows
\begin{equation}\label{eq:Dfwdcontr1}
\begin{aligned}
& \hat{v}_{i,h}^{n+1}=\hat{v}^n_{i,h}+\frac{\Delta t}{N}\sum_{j=1}^{N}H^n_{ij}\sum_{m=0}^M \hat{K}^n_{mh}(\hat{v}^n_{j,m}-\hat{v}^n_{i,m})+\Delta t \hat{u}_h^nQ_{i,h}^n,\\
& \hat{v}^n_{i,h}=\bar{\hat{v}}_{i,h},\\
& \hat{u}_h^n = -\frac{\Delta t}{\nu N} \sum_{i=1}^N(\hat{v}_{i,h}^{n+1}-\hat{v}_{d,h})Q_{i,h}^n,
\end{aligned}
\end{equation}
where $H^n_{ij}\equiv H(x^n_i,x^n_j)$ and $Q_{i,h}^n\equiv Q(\hat{v}^n_{i,h})$. Note that since the feedback control $\hat{u}_h^n$ in \eqref{eq:Dfwdcontr1} depends on the velocities at time $n+1$, the constrained interaction at time $n$ is implicitly defined. The feedback controlled system in the discretized form results
\begin{equation*}\label{eq:Dfwdcontr2}
\begin{aligned}
 \hat{v}^{n+1}_{i,h}&=\hat{v}^n_{i,h}+\frac{\Delta t}{N}\sum_{j=1}^{N}H^n_{ij}\sum_{m=0}^M \hat{K}^n_{mh}(\hat{v}^n_{j,m}-\hat{v}^n_{i,m})-\frac{\Delta t^2}{\nu N} \sum_{j=1}^N(\hat{v}_{j,h}^{n+1}-\hat{v}_{d,h})Q_{j,h}^nQ_{i,h}^n,\\ 
 \hat{v}^n_{i,h}&=\bar{\hat{v}}_{i,h}.
\end{aligned}
\end{equation*}
Again the action of the control  is substituted by an implicit term representing the relaxation toward the desired component of the velocity $\hat{v}_{d,h}$, and it can be inverted in a fully explicit system.

Considering the scaling for the regularization parameter $\nu=\kappa\Delta t $, the previous scheme is a consistent discretization of the following continuos system
\begin{equation}\begin{split}
\dfrac{d}{dt}\hat{v}_{i,h}(t)=&\dfrac{1}{N}\displaystyle\sum_{j=1}^NH(x_i,x_j)\sum_{m=0}^M \hat{K}_{mh}(t)(\hat{v}_{j,m}(t)-\hat{v}_{i,m}(t))\\
&+\dfrac{1}{\kappa N}\sum_{j=1}^N (\hat{v}_{d,h}-\hat{v}_{j,h}(t))Q(\hat{v}_{j,h}(t))Q(\hat{v}_{i,h}(t)).
\end{split}
\label{eq:CScontrolled}
\end{equation}
Now the control is explicitly embedded in the dynamic of the $h$th component of the gPC approximation as a feedback term, and the parameter $\kappa>0$ determines  its strength.

\subsection{Choice of the selective control}
For the specific choice of weight function $Q(\cdot)\equiv 1$ we refer in general to non selective control.  Note that in this case the action of the control is not strong enough in order to control the velocity of each agent, indeed in this case we are able only to control the mean velocity of the system. In  fact the control term is reduced to \begin{equation} \dfrac{1}{\kappa}(\hat{v}_{d,h}-\hat{\vm}_h), \end{equation} where $\hat{\vm}_h$ is the $h-$th coefficient of the expansion of $\vm$, that is
\begin{equation*}
\hat{\vm}_h = \dfrac{1}{N}\sum_{j=1}^N \hat{v}_{j,h}(t).
\end{equation*}
Then, only the projections of the mean velocity are steered toward the respective components of the target velocity, i.e. as soon as $\kappa\rightarrow 0$ it follows that $\hat{\vm}_h=\hat{v}_{d,h}$. Therefore, the choice of the selective function $Q(\cdot)$ is of paramount importance to ensure the action of the control on the single agent.

In principle one can address directly the control problem on the original system \eqref{eq:CS} as
\begin{equation}\begin{cases}\label{eq:CSconstrained}
\dot{x}_i(\theta,t)=v_i(\theta,t)\\
\dot{v}_i(\theta,t)=\dfrac{K(\theta,t)}{N}\displaystyle\sum_{j=1}^NH(x_i,x_j)(v_j(\theta,t)-v_i(\theta,t))+uQ(v_i(\theta,t)),
\end{cases}\end{equation}
where the control term $u$ is solution of 
\begin{equation}
u=\textnormal{arg}\min_{u}\left[\dfrac{1}{2}\int_0^T\dfrac{1}{N}\sum_{i=1}^N(v_i(\theta,t)-v_d)^2 dt+\dfrac{\nu}{2}\int_0^Tu(t)^2dt\right].
\end{equation}
Here $v_d\in \mathbb{R}^d$ is a target velocity and $\nu>0$ a regularization parameter. 
Similarly to previous subsection, through the approach presented in \cite{AHP,AP15,APZ}, we can derive the time-continuos MPC formulation which explicitly embed the control term in the dynamic as follows
\begin{equation}\begin{cases}\label{MPC_system}
\dot{x}_i(\theta,t)=&v_i(\theta,t)\\
\dot{v}_i(\theta,t)=&\dfrac{K(\theta,t)}{N}\displaystyle\sum_{j=1}^NH(x_i,x_j)(v_j(\theta,t)-v_i(\theta,t))\\
&\quad+\dfrac{1}{\kappa N}\displaystyle\sum_{j=1}^N(v_d-v_j(\theta,t))Q(v_j(t,\theta))Q(v_i(t,\theta)).
\end{cases}\end{equation}
Now the gPC approximation of \eqref{MPC_system} can be obtained as in Section \ref{sec:gPC} and leads to the set of ODEs
\begin{equation}\begin{split}\label{MPC_Ncontrolled}
\dfrac{d}{dt}\hat{v}_{i,h}(t)=&\dfrac{1}{N}\displaystyle\sum_{j=1}^N H(x_i,x_j)\sum_{m=0}^M \hat{K}_{mh}(t)(\hat{v}_{j,m}(t)-\hat{v}_{i,m}(t))\\
&+\dfrac{1}{\kappa N}\sum_{j=1}^N R_{h}(v^M_i,v^M_j),
\end{split}
\end{equation}
where 
\begin{equation}\label{eq:Rmh}
R_{h}(v^M_i,v^M_j) =\dfrac{1}{\|\Phi_h \|^2}\mathbb{E}_{\theta}\left[(v_d-v_j^M)Q\left(v_i^M(\theta,t)\right)Q\left(v_j^M(\theta,t)\right)\Phi_h(\theta)\right].
\end{equation} 
In general systems \eqref{MPC_Ncontrolled} and  \eqref{eq:CScontrolled}, without further assumptions on the selective function $Q(\cdot)$, are not equivalent. In addition to the non selective case, there exist at least one choice of selective control that makes the two approaches totally interchangeable. In fact, taking 
\begin{equation}\label{eq:Q_linear}
Q(v_i)=\dfrac{v_d-v_i}{\sqrt{\frac{1}{N}\sum_{j=1}^N(v_d-v_j)^2}},
\end{equation}
we have $Q(\cdot)\in [-1,1]$ and the control term in (\ref{MPC_Ncontrolled}) takes the following form
\begin{equation}\begin{split}\label{eq:feedback1}
\dfrac{1}{\kappa N}\sum_{j=1}^N R_{h}=\dfrac{1}{\kappa ||\Phi_h ||^2}\mathbb{E}_{\theta}\left[(v_d-v_i^M(\theta,t))\Phi_h(\theta)\right] 
 = \dfrac{1}{\kappa}\left(\hat{v}_{d,h}-\hat{v}_{i,h}\right).
\end{split}\end{equation}
Similarly the control term in (\ref{eq:CScontrolled}) reduces to
\begin{equation}\begin{split}\label{eq:feedback2}
\dfrac{1}{\kappa N}\sum_{j=1}^N (\hat{v}_{d,h}-\hat{v}_{j,h}(t))Q(\hat{v}_{j,h}(t))Q(\hat{v}_{i,h}(t))=\dfrac{1}{\kappa}\left(\hat{v}_{d,h}-\hat{v}_{i,h}\right),
\end{split}\end{equation}
and therefore system \eqref{MPC_Ncontrolled} coincides with \eqref{eq:CScontrolled}. Note that as $\kappa \to 0$ both systems are driven towards the controlled state $\hat{v}_{i,h}=\hat{v}_{d,h}$ which implies a strong control over each single agent.

\begin{figure}[t]
\begin{center}
\begin{tikzpicture}
[node distance=.5cm, start chain=going below,]
\node [punktchain] (model)  {\textnormal{Model}};
\node[below=2cm of model](dummy){};
\begin{scope}[start branch=venstre,every join/.style={->, thick, shorten <=1pt}, ]
\node [punktchain,left=of dummy] (control_1)  {\textnormal{Control Problem}};
\node [punktchain,join,on chain=going below] (MPC_1)  {\textnormal{MPC}};
\node [punktchain,join, on chain=going below] (gPC_1)  {\textnormal{gPC}};
\end{scope}
\begin{scope}[start branch=hoejre,every join/.style={->, thick, shorten <=1pt}, ]
\node[punktchain,right=of dummy]  (gPC_2)  {\textnormal{gPC}};
\node [punktchain,join, on chain=going below] (control_2)  {\textnormal{Control Problem}};
\node [punktchain,join, on chain=going below] (MPC_2)  {\textnormal{MPC}};
\end{scope}
\node[below=2cm of gPC_1](dummy_2){};
\node [punktchain,right=of dummy_2] (constrained control)  {\textnormal{Constrained gPC System}};
\draw[|-,-|,->, thick,](model.south) |-+(0,-1em)-| (control_1.north);
\draw[|-,-|,->, thick,] (model.south) |-+(0,-1em)-| (gPC_2.north);
\draw[thick,->](gPC_1.south) |-+(0,-3em)-| (constrained control.north);
\draw[thick,->](MPC_2.south) |-+(0,-3em)-| (constrained control.north);
\end{tikzpicture}
\caption{Sketch of the two numerical approaches to solve the control problem with uncertainty, combining MPC and gPC. In both cases, of non selective control, i.e. $Q(\cdot)\equiv1$, and of selective control with $Q(\cdot)$ defined in \eqref{eq:Q_linear} the two approaches are equivalent.}
\label{fig:scheme}
\end{center}
\end{figure}
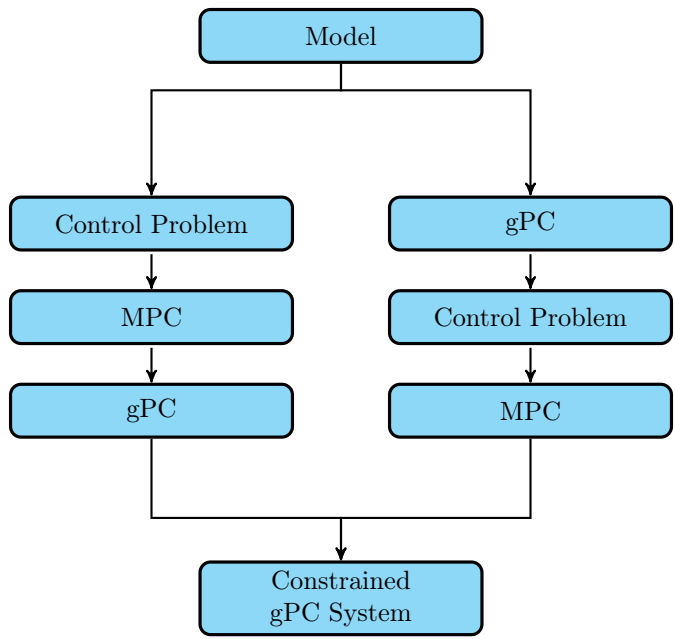
In Figure \ref{fig:scheme} we summarize the two equivalent approaches. In the case of non selective control and of selective function given by \eqref{eq:Q_linear} the constrained gPC system can be obtained from our initial unconstrained model \eqref{eq:CS} through two different but equivalent methods. The first approximates the solution of the Cucker-Smale type model via the gPC projection and then introduces a control on the coefficients of the decomposition through a MPC approach in order to steer each component to $(\hat{v}_{d,0},\hat{v}_{d,1},\dots,\hat{v}_{d,M})$. Whereas the second method considers a constrained Cucker-Smale problem \eqref{eq:CSconstrained}, introduces its continuous MPC approximation and then computes the gPC expansion of the resulting system of constrained differential equations.
\begin{remark}
We remark that the choice of $Q(\cdot)$ stated in \eqref{eq:Q_linear}, for which the two approaches sketched in Figure \ref{fig:scheme} are identical, is equivalent to consider the constrained dynamic \eqref{eq:CSconstrained}, modified as follows
\begin{equation}\begin{cases}\label{eq:CSconstrainedUi}
\dot{x}_i(\theta,t)=v_i(\theta,t)\\
\dot{v}_i(\theta,t)=\dfrac{K(\theta,t)}{N}\displaystyle\sum_{j=1}^NH(x_i,x_j)(v_j(\theta,t)-v_i(\theta,t))+u_i,
\end{cases}\end{equation}
where the control term, $u_i$ for each agent $i=1,\ldots,N$, is given by the minimization of the following functional
\begin{align}\label{eq:funcJ}
J(v_1,\ldots,v_N;u_1,\ldots,u_N)=\frac{1}{2}\int_0^T\frac{1}{N}\sum_{i=1}^N\left[(v_i(\theta,t)-v_d)^2+\frac{\nu}{2}u_i(t)^2\right]dt.
\end{align}
Since the functional is strictly convex, applying the (MPC) procedure to \eqref{eq:CSconstrainedUi}-\eqref{eq:funcJ}, we obtain as first order approximation for the solution of the optimal control problem the feedback control term\begin{align}
u_i = \frac{1}{\kappa}(v_d-v_i),\qquad i=1,\ldots,N.
\end{align}
Thus the same considerations on the equivalence of the approaches hold.
\end{remark} 

\section{Numerical tests}\label{sec:Numerics}
\label{sec:Numeric}
We present some numerical experiments of the behavior of the flocking model in the case of a Hermite polynomial chaos expansion. This choice corresponds to the assumption of a normal distribution for the stochastic parameter in the Cucker-Smale type equation \eqref{eq:CS} and in its constrained behavior \eqref{eq:CScontrolled}. Numerical results show that the introduced selective control with the weight function \eqref{eq:Q_linear} is capable to drive the velocity to a desired state even in case of a dynamic dependent by a normally distributed random input, with fixed or time-dependent variance. In the uniform interaction case, since the effect of agents' positions do not influence the alignment we report only the results of the agents' velocities.

\begin{figure}[t]
\begin{center}
\hskip -.4cm
\!\!
\includegraphics[scale=0.317]{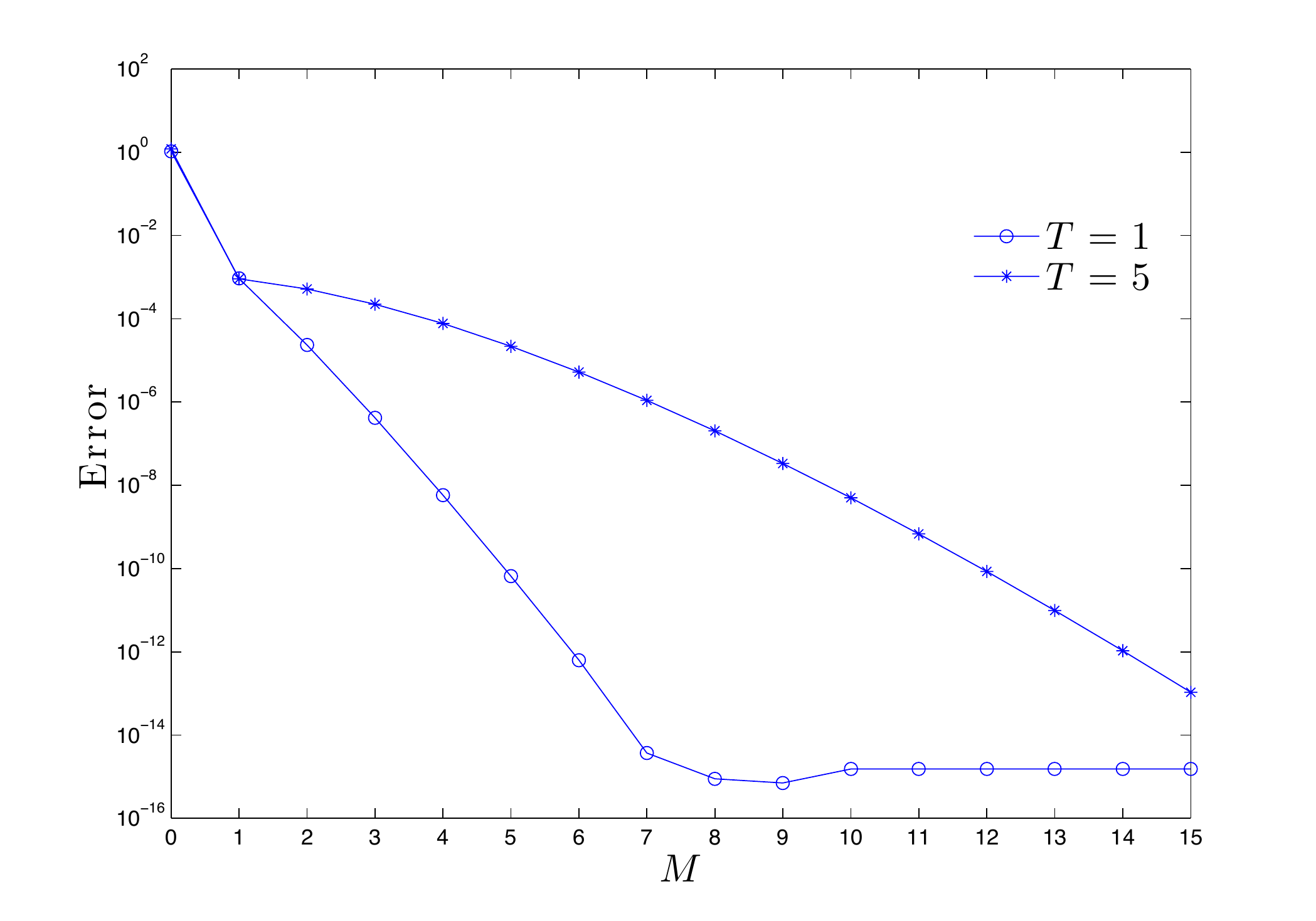}
\hskip -0.75cm
\includegraphics[scale=0.317]{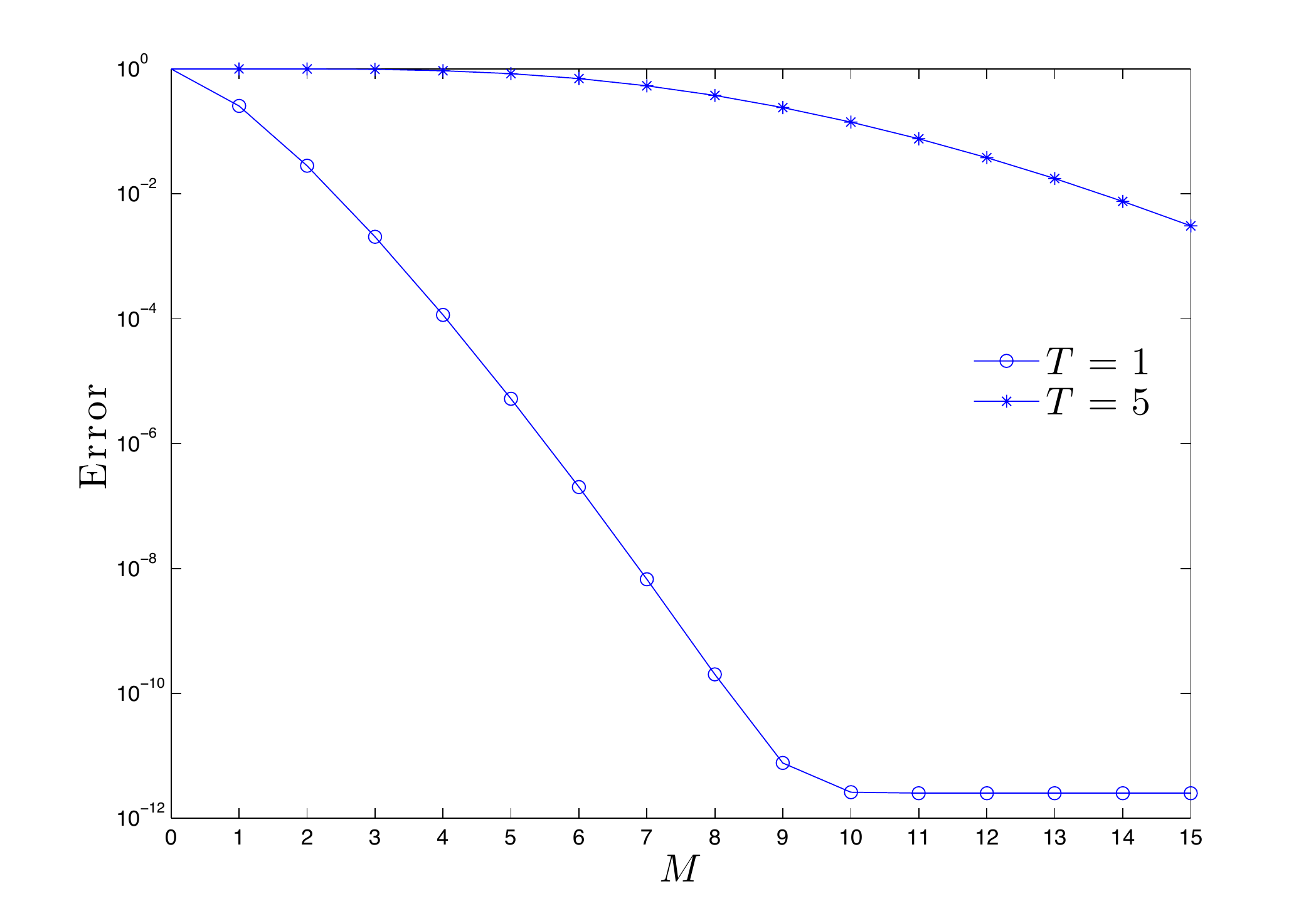}\!\!\!\!\!\!
\caption{Error convergence for increasing number of polynomials in the gPC decomposition approximation. Left: convergence of the mean error at two fixed times $T=1$ and $T=5$. Right: convergence of the variance error. In both cases we considered a random time-independent scattering $K(\theta,t)=\theta$, where the random variable $\theta$ is normally distributed $N(2,1/2)$. The system of ODEs is solved through a $4$th order Runge-Kutta with $\Delta t=10^{-5}$.}
\label{Fig:convergence}
\end{center}
\end{figure}

\begin{figure}[t]
\begin{center}
\includegraphics[scale=0.32]{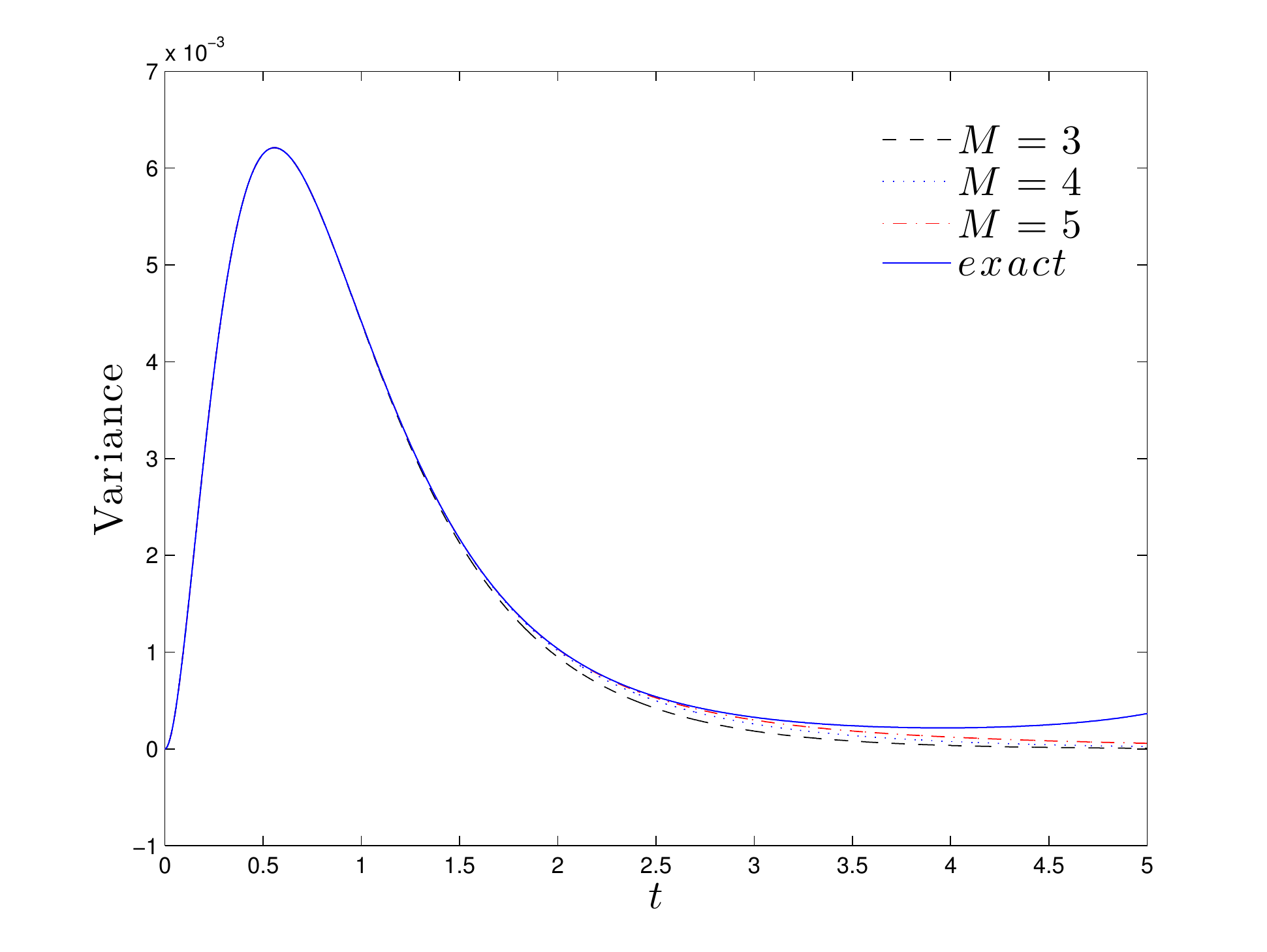}\,\,
\hskip -0.72 cm
\includegraphics[scale=0.32]{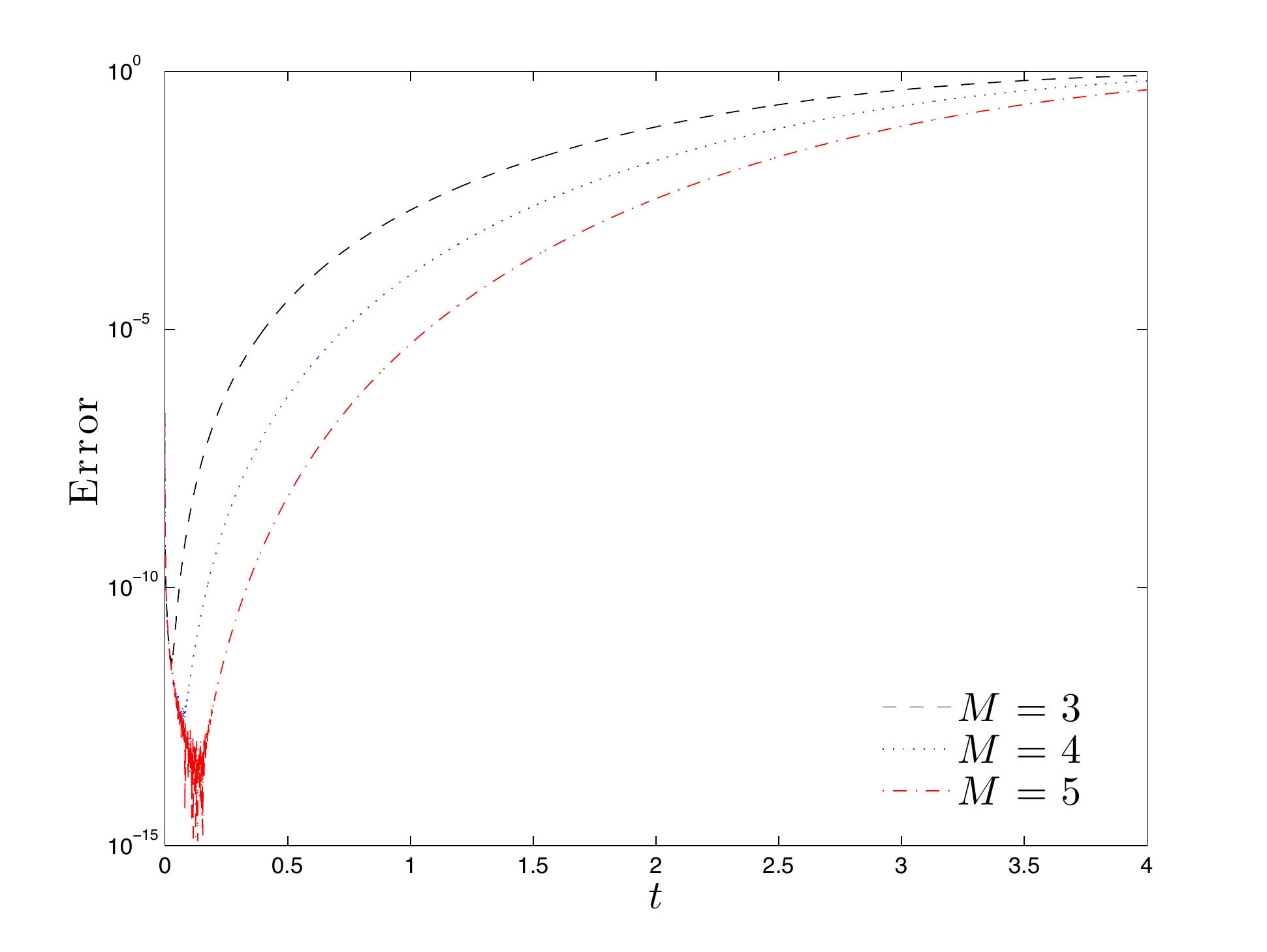}
\caption{Evolution of the variance-error $E_v(t)$ defined in equation \eqref{eq:errors} for the gPC decomposition for the unconstrained model \eqref{eq:CShomo} with $K(\theta,t)=\theta\sim\mathcal{N}(2,1/2)$ over the time interval $[0,T]$ with $T=5$ and time step $\Delta t = 10^{-5}$. }
\label{Fig:error_variance}
\end{center}
\end{figure}

\begin{figure}[h]
\begin{center}
\includegraphics[scale=0.32]{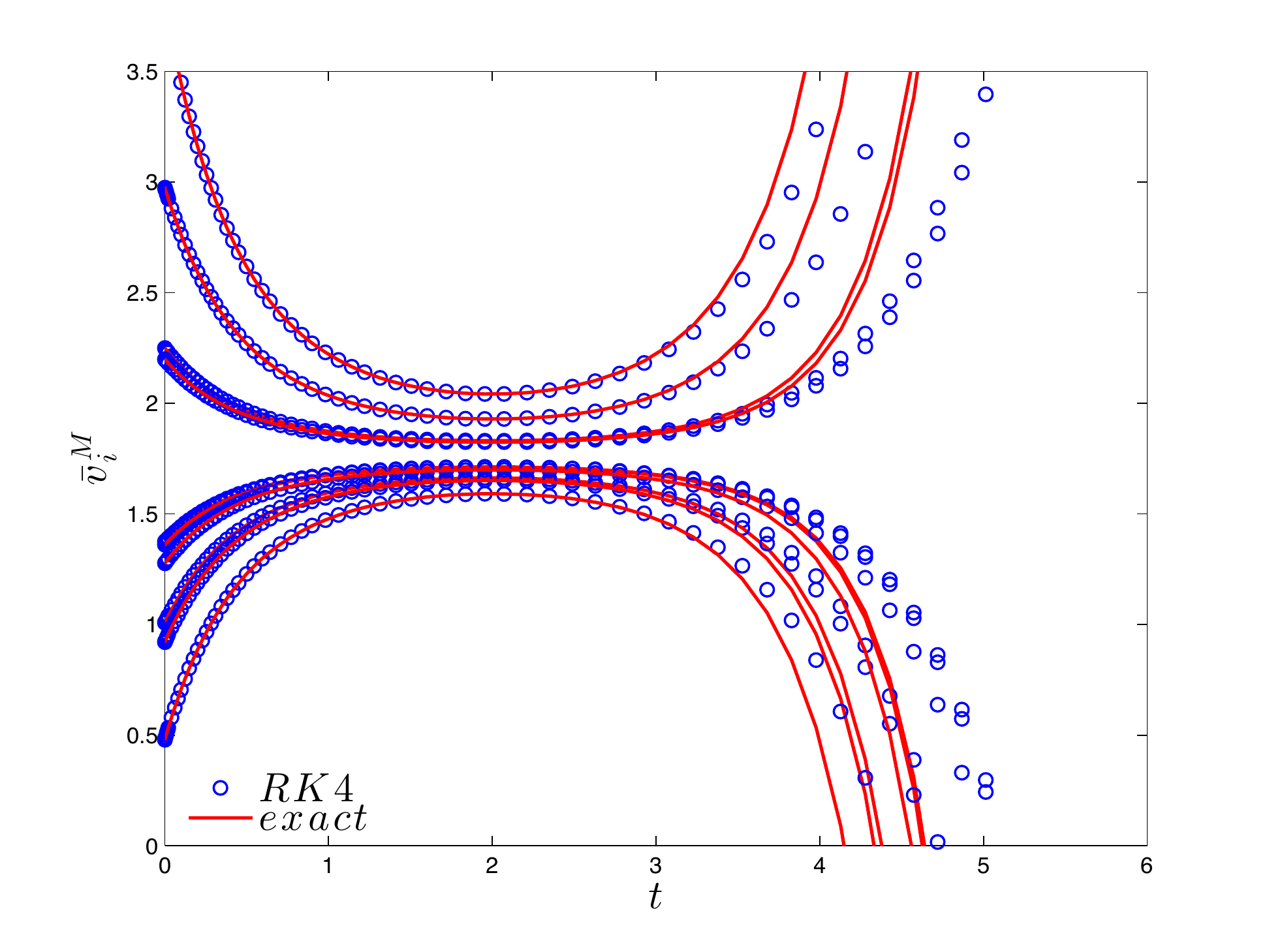}
\hskip -0.6 cm
\includegraphics[scale=0.32]{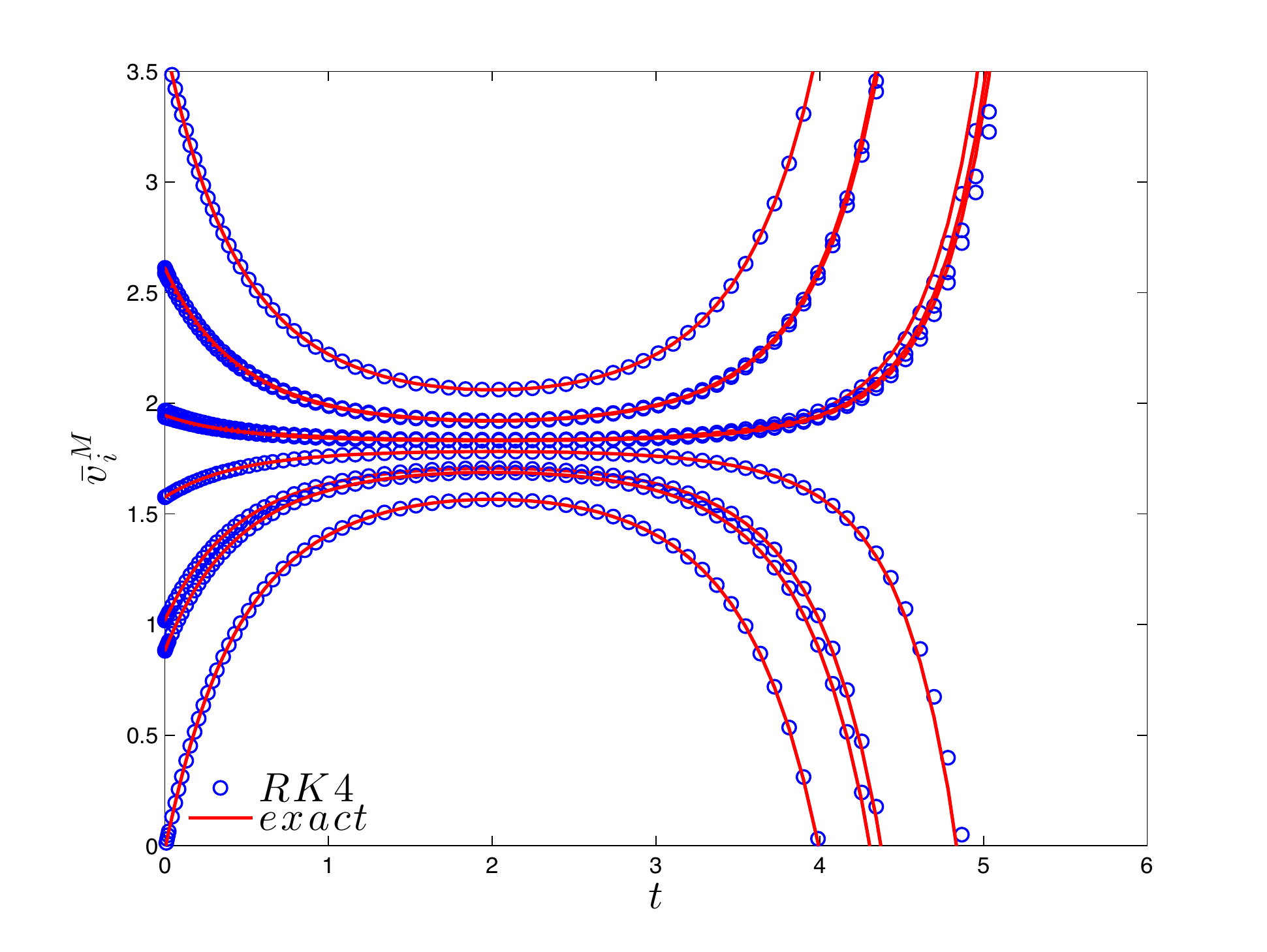}
\caption{Left: $6$th order Hermite gPC decomposition solved through a $4$th order Runge-Kutta. Right:  $10$th order Hermite gPC decomposition solved through a $4$th order Runge-Kutta. In both cases  the final time considered is $T= 6$, with time step $\Delta t = 10^{-2}$. }
\label{Fig:errors}
\end{center}
\end{figure}

\subsection{Unconstrained case}
In Figures \ref{Fig:convergence} and \ref{Fig:error_variance} we present numerical results for the convergence of the error using the gPC scheme described in equation \eqref{eq:decomposition} for $H\equiv 1$ and solved through a $4$th order Runge-Kutta method. In particular Figure \ref{Fig:convergence} shows the behavior of the error with respect to increasing terms of the gPC decomposition. Here we considered the average in time of the error for the mean and the variance at time $t>0$ in the $L^1$ norm 
\begin{equation}\label{eq:errors}
E_{\bar{v}}(t) = \dfrac{1}{N}\sum_{i=1}^N\left | \dfrac{\bar{v}_i(t)-\bar{v}^M_i(t)}{\bar{v}_i(t)}\right | \qquad E_{\bar{\sigma}^2}(t)=\dfrac{1}{N}\sum_{i=1}^N \left | \dfrac{\bar{\sigma}^2_i(t)-\bar{\sigma}^{2,M}_i(t)}{\bar{\sigma}^2_i(t)} \right |,
\end{equation}
where
\begin{equation}
\bar{\sigma}^2_i(t)=\mathbb{E}_{\theta}\left[\left(v_i(\theta,t)-\bar{v}_i(t)\right)^2\right]
\end{equation}
with $v_i(\theta,t)$ and $\bar{v}_i(t)$ defined in \eqref{eq:exact_homo} and \eqref{eq:expected_velocity}. Observe that if the scattering rate $K(\theta,t)$ is of the from described in \eqref{eq:K_normal} with $h(\cdot)\equiv 1$ and $k(\theta)\sim\mathcal{N}(\mu,\sigma^2)$ than, in addition to the explicit evolution for the expected velocity as in \eqref{eq:expectation}, we can obtain the exact version for the evolution of the variance of the $i$th agent
\begin{equation}
\bar{\sigma}^2_i(t)=(v_i(0)-\vm )^2 \left(\exp\{{-2\mu t+2\sigma^2 t^2}\}-\exp\{{-2\mu t+\sigma^2 t^2}\}\right).
\end{equation}
In \eqref{eq:errors} we indicated with $\bar{\sigma}^{2,M}_i(t)$ the approximated variance
\begin{equation}
\bar{\sigma}^{2,M}_i(t)=\sum_{h=0}^M \hat{v}^2_{i,h}(t)\mathbb{E}_{\theta}[\Phi_h(\theta)^2]-\hat{v}^2_{i,0}(t).
\end{equation}
It is easily seen how the error decays spectrally for increasing value of $M$, however the method is not capable to go above a certain accuracy and therefore for large $M$ a threshold effect is observed. This can be explained by the large integration interval we have considered in the numerical computation, and by the well-known loss of accuracy of gPC for large times \cite{Vos}. In the case of the error of the variance, Figure \ref{Fig:error_variance}, the gPC approximation exhibits a  slower convergence with respect to the convergence of the mean. Next in Figure \ref{Fig:errors} we see how for large times the solution of the differential equation \eqref{eq:CShomo} diverges and the numerical approximation is capable to describe accurately its behavior only through an increasing number of Hermite polynomials. 

\subsection{Constrained uniform interaction case}
In  Figure \ref{Fig:ExactControlled} we show different scenarios for the uniform interaction dynamic with constraints.  
In the first row we represents the solution for $N=10$ agents, whose dynamic is described by equation \eqref{eq:CScontrolled} with $v_d=1$, different values of $\kappa$ originate different controls on the average of the system, which however do not prevent the system to diverge.
In the second row we show the action of selective control \eqref{eq:Q_linear}.
It is evident that, with this choice, we are able to control the system also in the case with higher variance. 

Observe that the numerical results are coherent with the explicit solution of the controlled equation. Let us consider the time-independent scattering rate$K(\theta,t)=\theta\sim\mathcal{N}(\mu,\sigma^2)$, then from the equation
\begin{equation}\label{eq:controlled_differential}
\dfrac{d}{dt}v_i(\theta,t)=\theta(\vm -v_i(\theta,t))+\dfrac{1}{\kappa}(v_d-v_i(\theta,t))
\end{equation}
we can compute the exact solution given $v_i(\theta,0)=v_i(0)$
\begin{equation}\begin{split}\label{eq:eq_controlled_exact}
v_i(\theta,t)=&\dfrac{\kappa\vm \theta  + v_d}{\kappa \theta+1}+\left(v_i(0)-\dfrac{\kappa\vm  \theta + v_d}{\kappa \theta+1}\right)\exp{\Big\{-\left(\theta+\dfrac{1}{\kappa}\right)t\Big\}}.
\end{split}\end{equation}
The asymptotic behavior of the expected value of \eqref{eq:eq_controlled_exact} can be studied similarly to what we did in Section \ref{uniform}. In other words in order to prevent the divergence of the leading term of the controlled expected exact solution we might study 
\begin{equation}\label{eq:esponente_controlled}
\exp{\Big\{-\left(\mu+\dfrac{1}{\kappa}\right)t+\dfrac{\sigma^2 t^2}{2}\Big\}},
\end{equation}
which diverge if 
\begin{equation}
t>\dfrac{2}{\sigma^2}\left(\mu+\dfrac{1}{\kappa}\right).
\end{equation}
Then for each fixed time we could select a regularization parameter $\kappa>0$ so as to avoid the divergence of \eqref{eq:eq_controlled_exact}. Moreover we can observe that in the limit $\kappa\rightarrow 0$ the introduced selective control is capable to correctly drive the system \eqref{eq:controlled_differential} for each $t>0$.

\begin{figure}[htb]
\begin{center}
\includegraphics[scale=0.32]{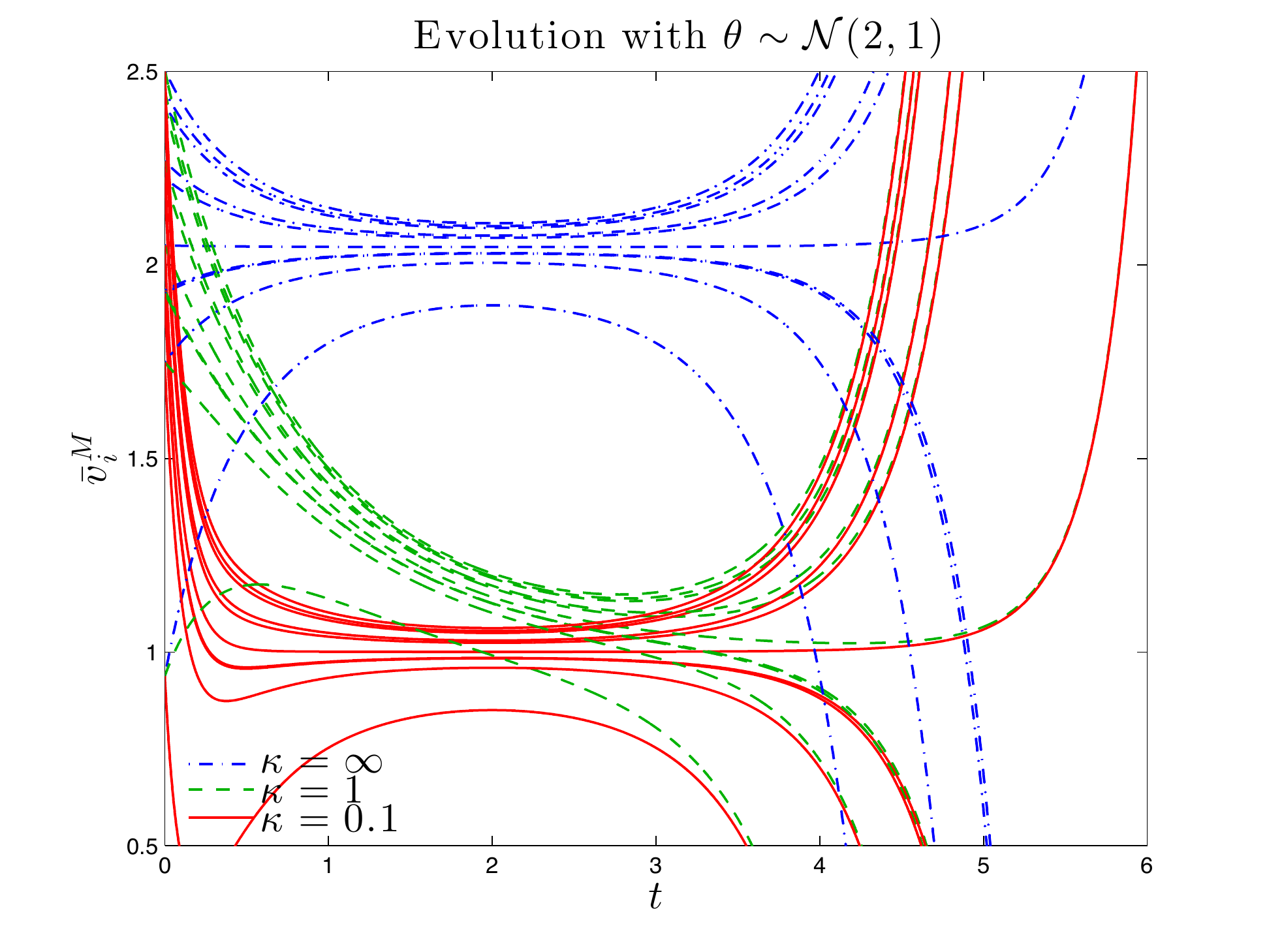}\,\,
\hskip -0.8 cm
\includegraphics[scale=0.32]{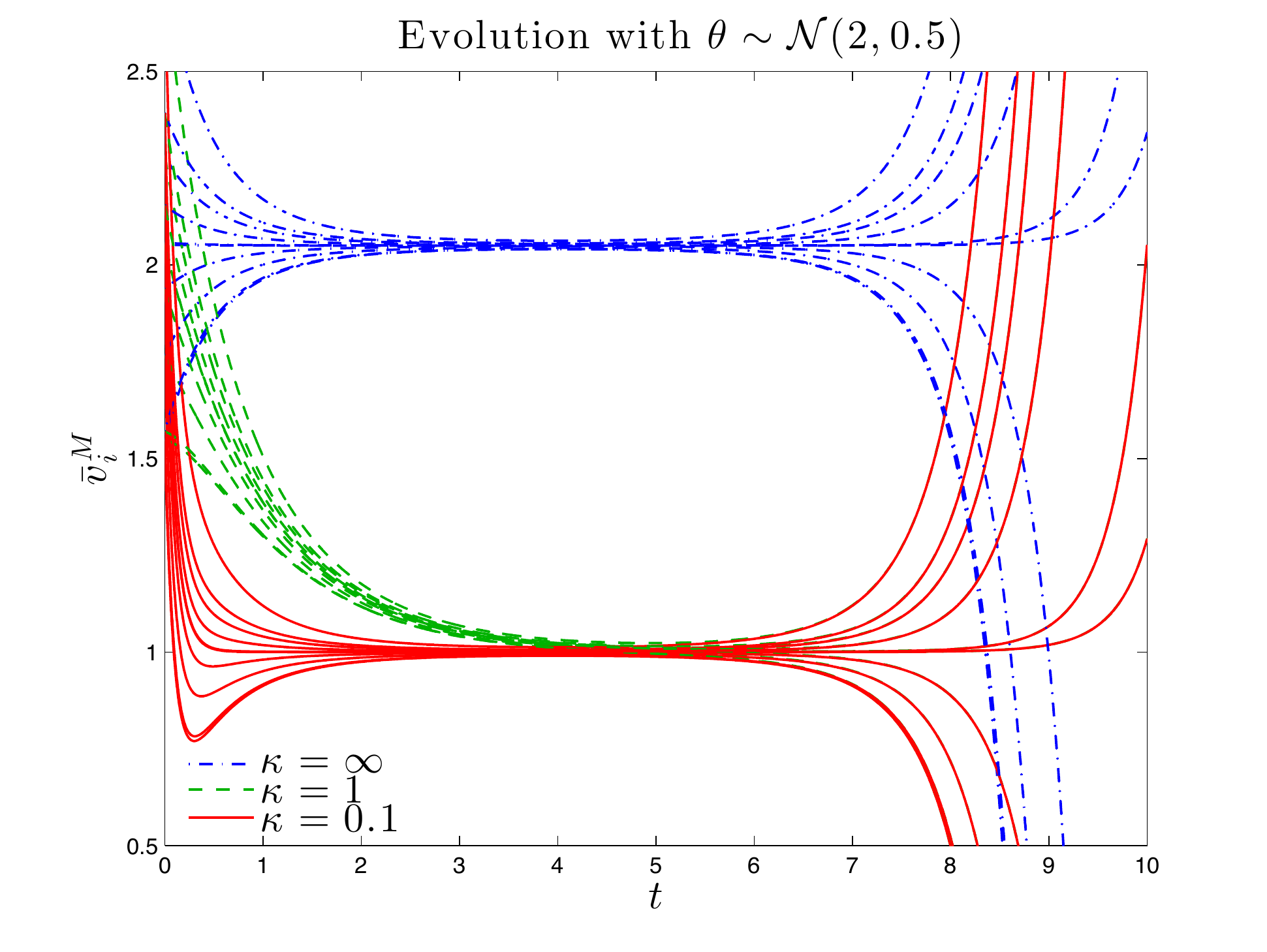}
\\
\includegraphics[scale=0.32]{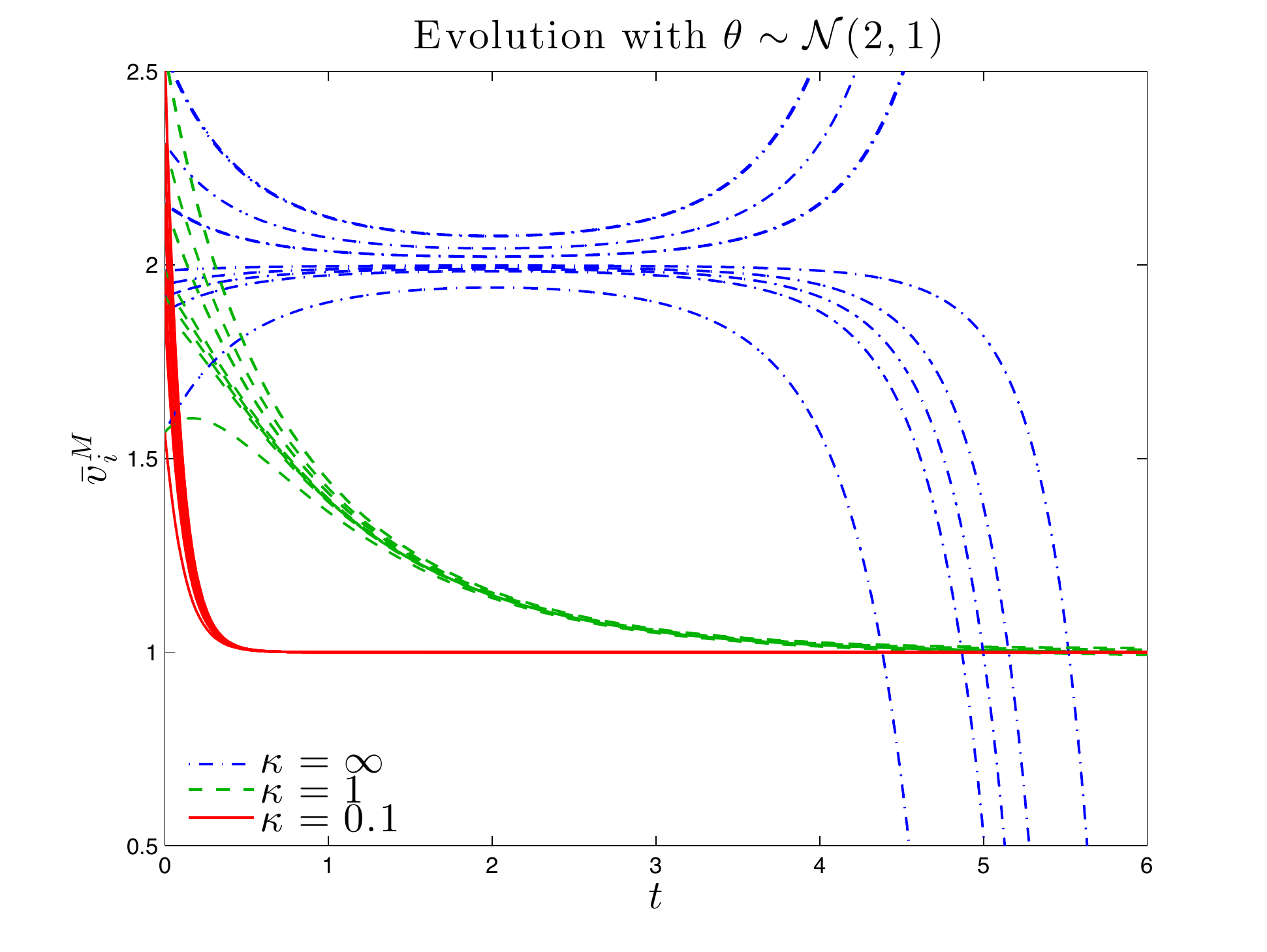}\,\,
\hskip -0.8 cm
\includegraphics[scale=0.32]{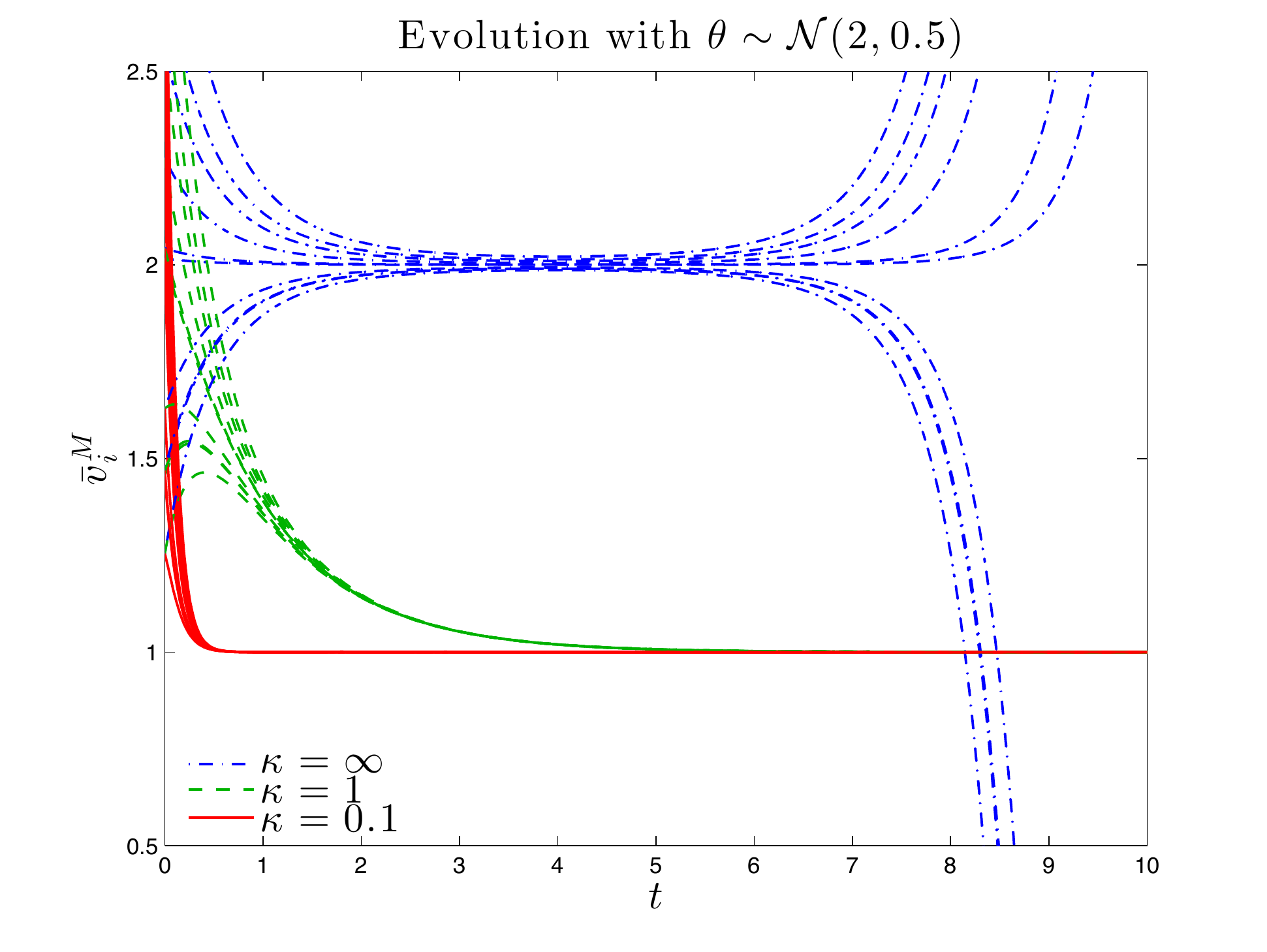}
\caption{Evolution of the uniform interaction alignment model \eqref{eq:CScontrolled} with $N=10$ agents, at $t=0$ distributed around $\vm =2$ with unitary variance, depending on a normal random parameter. Left column: $\theta\sim\mathcal{N}(2,1)$. Right column $\theta\sim\mathcal{N}(2,0.5)$. The control term shows its ability to steer the system towards desired velocity $v_d=1$, with different intensities $\kappa=1$ and $\kappa =0.1$, when $\kappa=\infty$ the control has no influence. First row shows the action of the control acting just on the average velocity, $Q\equiv 1$. Second row shows the action of selective control with $Q(\cdot)$ as in \eqref{eq:Q_linear}.  
}
\label{Fig:ExactControlled}
\end{center}
\end{figure}

In Figure \ref{Fig:variance_timedep} we consider the system with random time-dependent scattering rate $\theta\sim\mathcal{N}(\mu,\sigma^2(t))$. The dynamic shows how, for the choice of time dependent variance described in Remark \ref{rem:variance}, that is $\sigma(t)=1/s^\alpha$ with $\alpha=1/2$, the convergence depends from the mean value of the random input. In particular numerical experiments highlight the threshold effect for $\mu=2$ which we derived in Section \ref{sec:CS}. In the second figure we show that the action of the selective control \eqref{eq:Q_linear}, with desired velocity $v_d=\vm $, is capable to stabilize the system and drive the velocities towards the desired state.

\begin{figure}[]
\begin{center}
\includegraphics[scale=0.32]{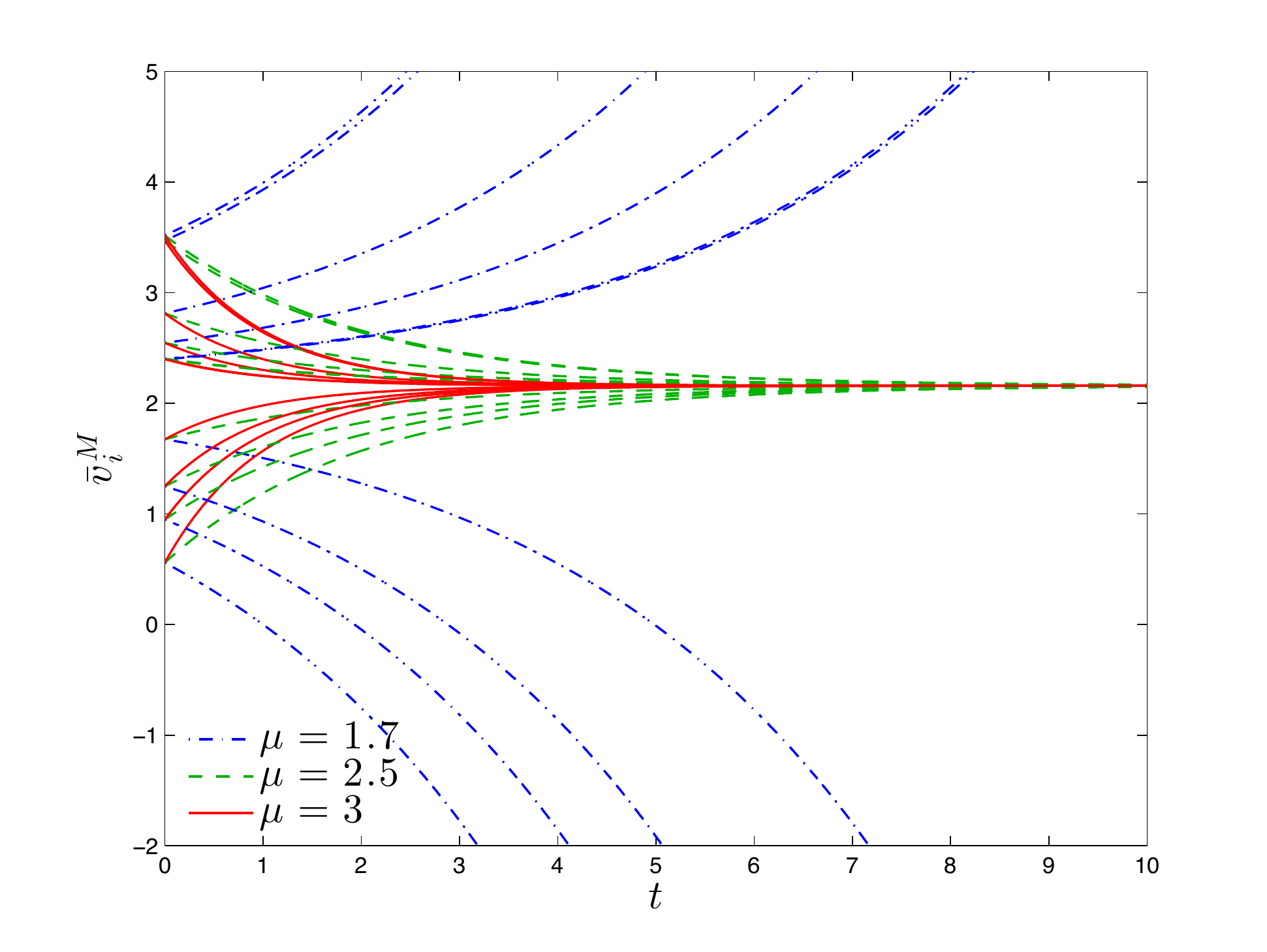}
\hskip -0.7cm
\includegraphics[scale=0.32]{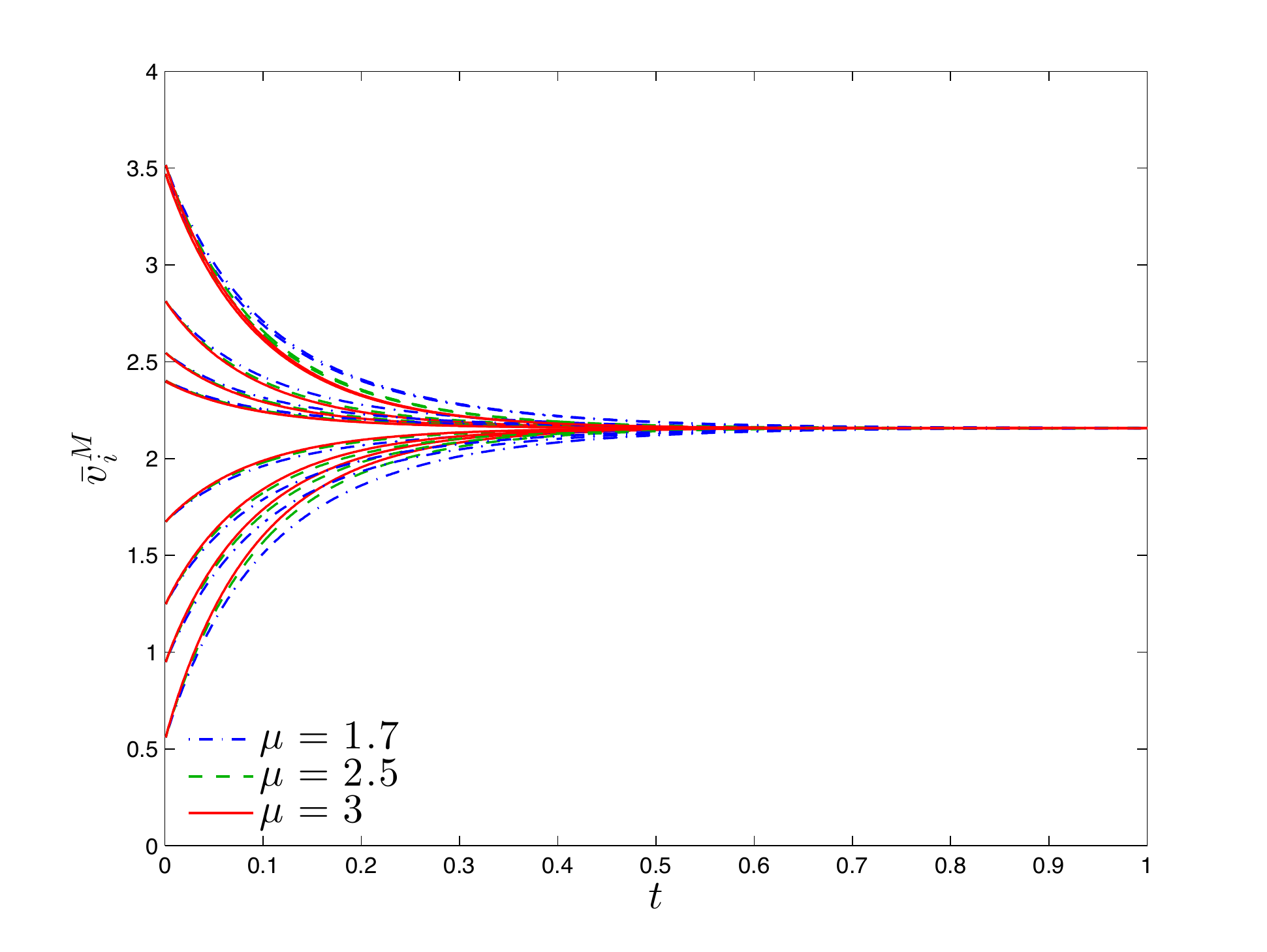}
\caption{Solution of the uniform interaction case with time dependent random parameter $\theta$ distributed accordingly to a normal distribution $\mathcal{N}(\mu,\sigma^2(t))$, with a time-dependent standard deviation $\sigma(t)=1/t^\alpha$, and $\alpha=1/2$. Left: we see the threshold for different values of $\mu$, i.e. for $\mu<2$ the system diverges. Right: solution of the constrained model with $\kappa=0.1$, observe that we are able to steer the system to the desired velocity $v_d=\vm $, i.e. the initial mean velocity of the system, using the selective control described in \eqref{eq:Q_linear}.}
\label{Fig:variance_timedep}
\end{center}
\end{figure}

\subsection{Constrained space dependent case}
Next let us consider the full space non homogeneous constrained problem \eqref{eq:CS} with the interaction function defined in \eqref{eq:interaction_space}. In this case we assume that $K(\theta)=\theta$ with $\theta\sim\mathcal{N}(\mu,\sigma^2)$. In Figure \ref{Fig:separation} and \ref{Fig:flocking} we consider a system of $N=100$ agents with Gaussian initial position with zero mean and with variance $2$ and Gaussian initial velocities clustered around $\pm 5$ with mean $1/10$ . The numerical results for \eqref{eq:decomposition} have been performed through a $10$th order gPC expansion. The dynamic has been observed in the time interval $[0,5]$ with $\Delta t=10^{-2}$.  In Figure \ref{Fig:flocking} we see how the selective control is capable to drive the velocity of each agent to the desired state $v_d$. In fact in case of no control, see Figure \ref{Fig:separation}, we have that the velocities of the system naturally diverges. 

\begin{figure}[tbh]
\begin{center}
\subfigure[t=0]{\includegraphics[scale=0.32]{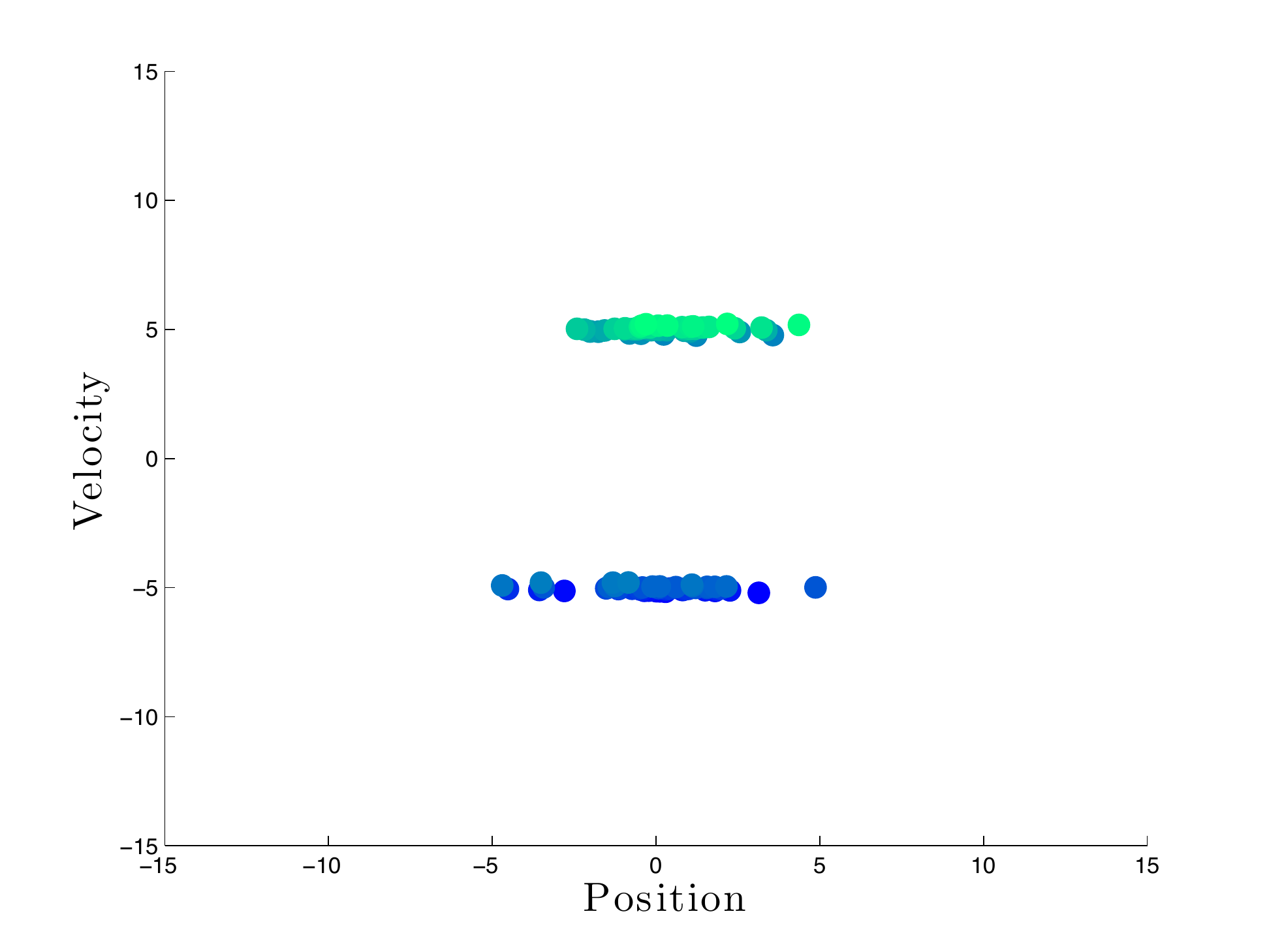}}
\hskip -0.7cm
\subfigure[t=1]{\includegraphics[scale=0.32]{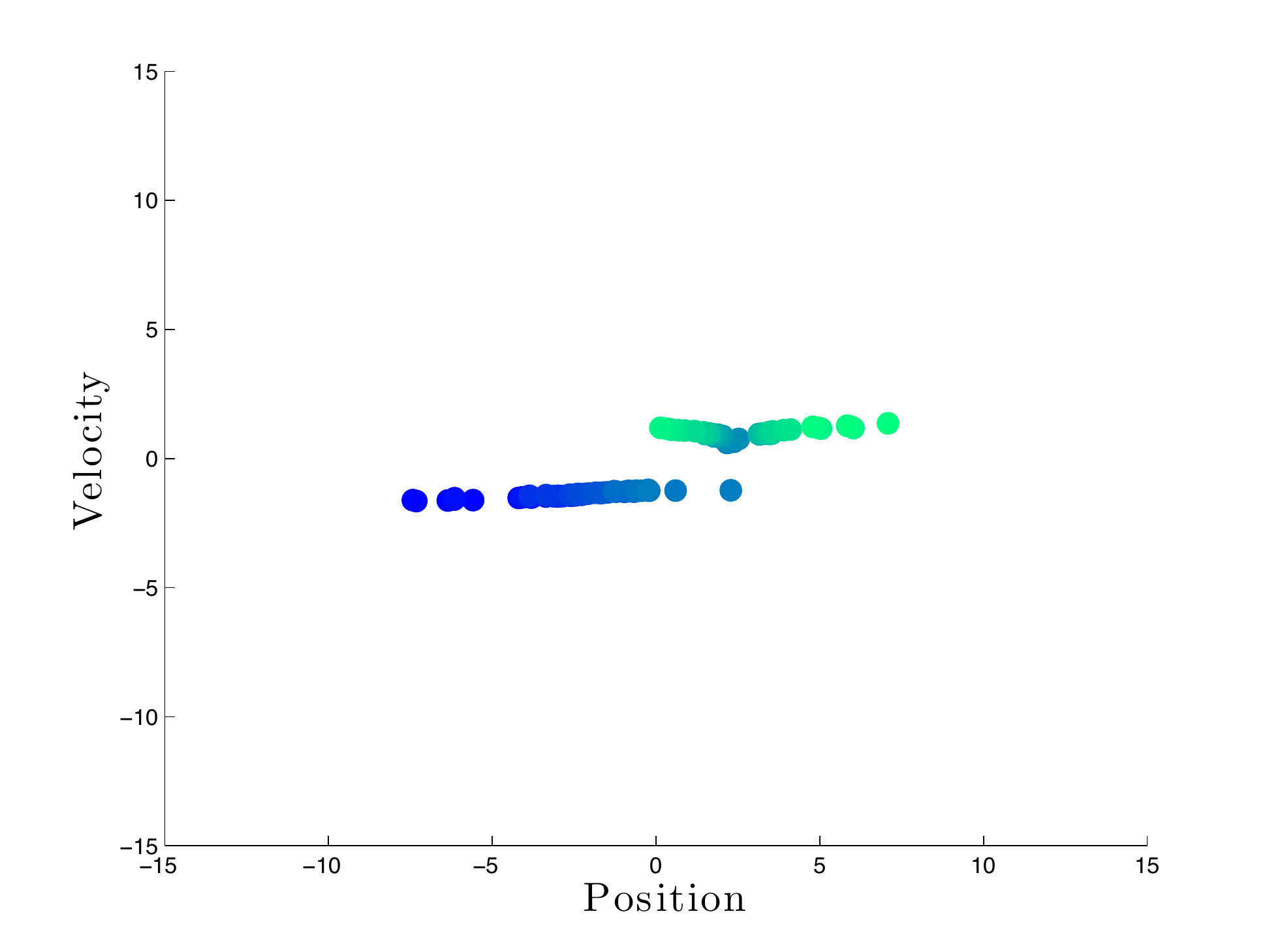}}\\
\subfigure[t=2]{\includegraphics[scale=0.32]{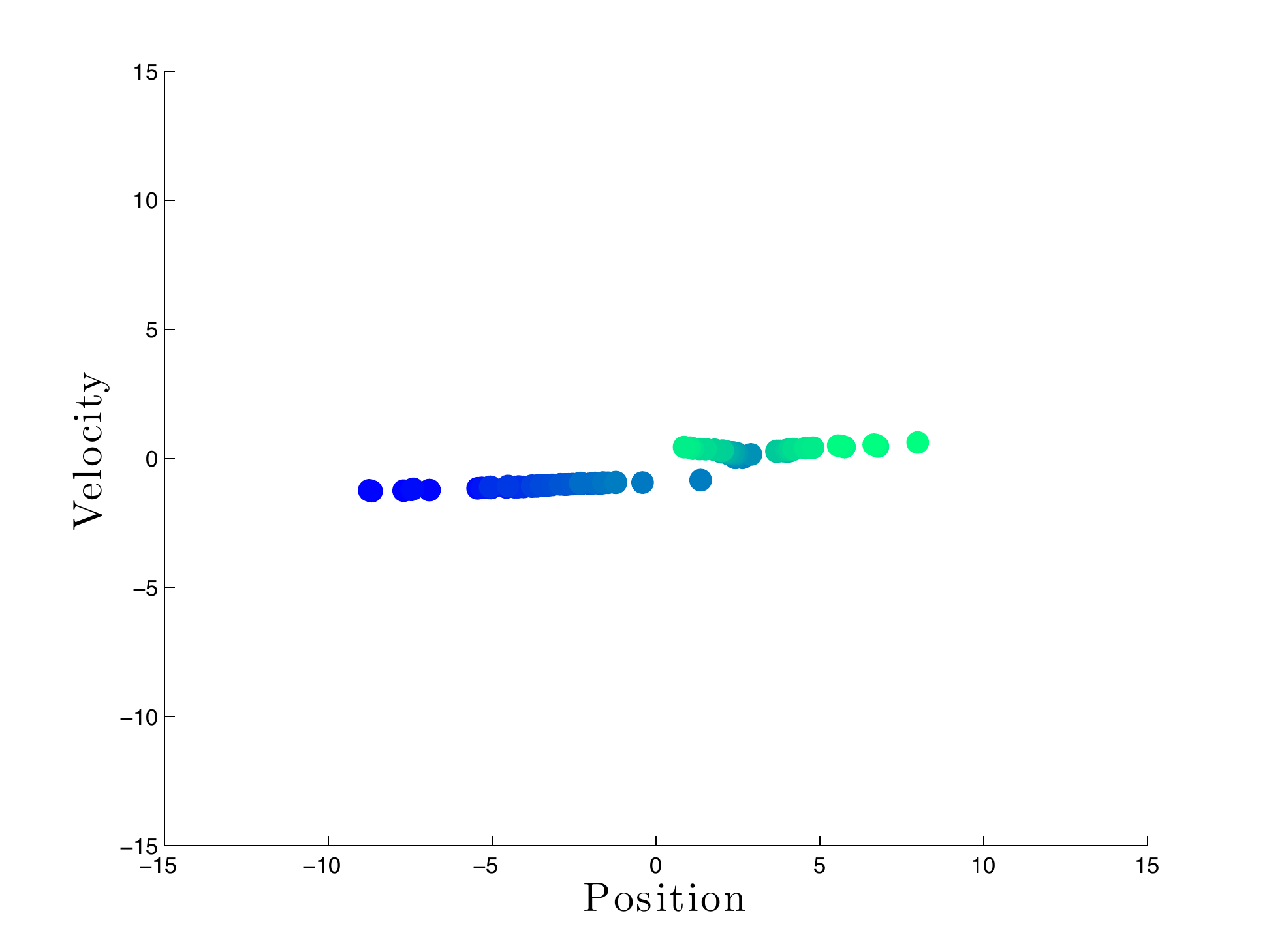}}
\hskip -0.7cm
\subfigure[t=3]{\includegraphics[scale=0.32]{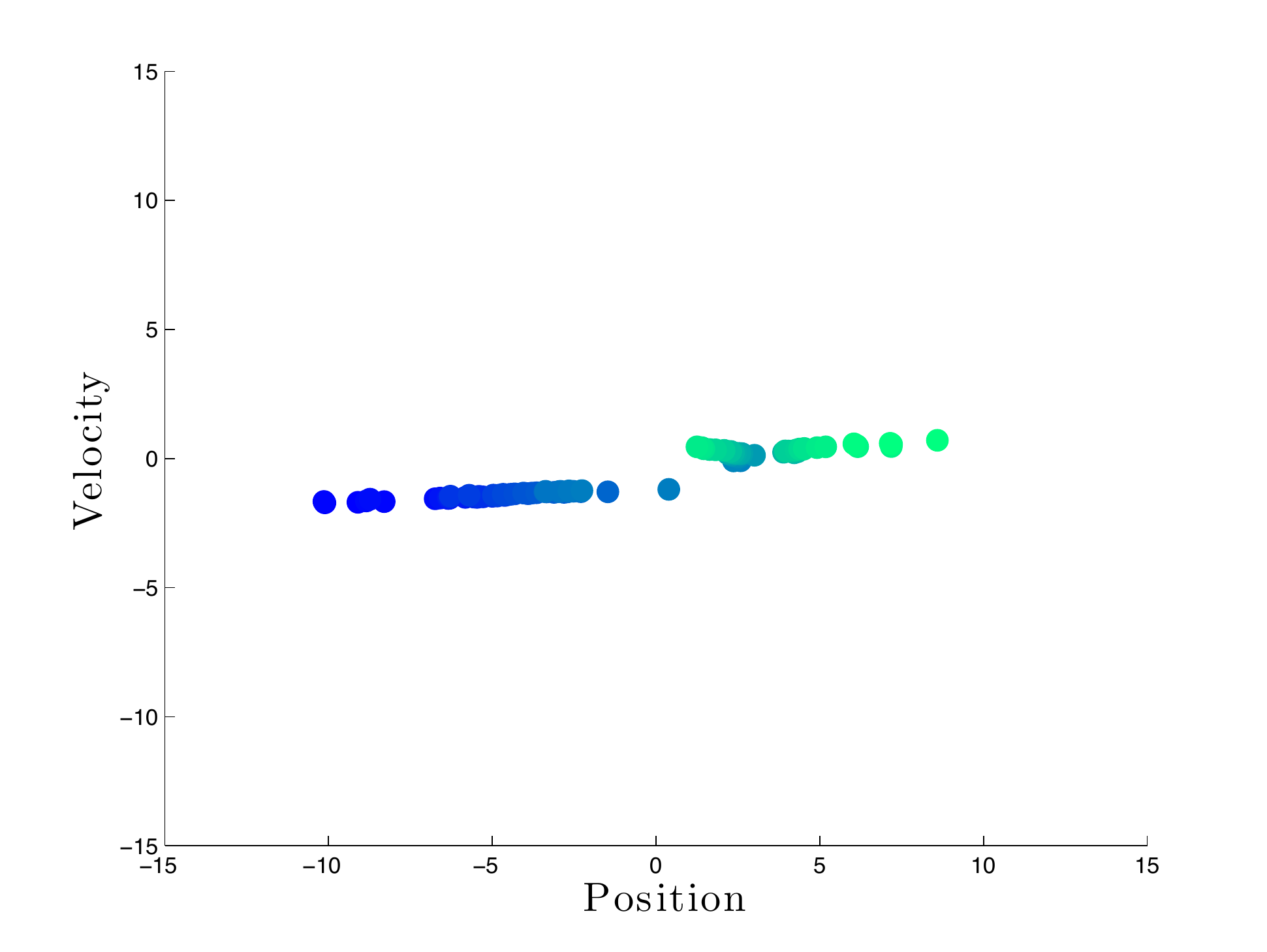}}\\
\subfigure[t=4]{\includegraphics[scale=0.32]{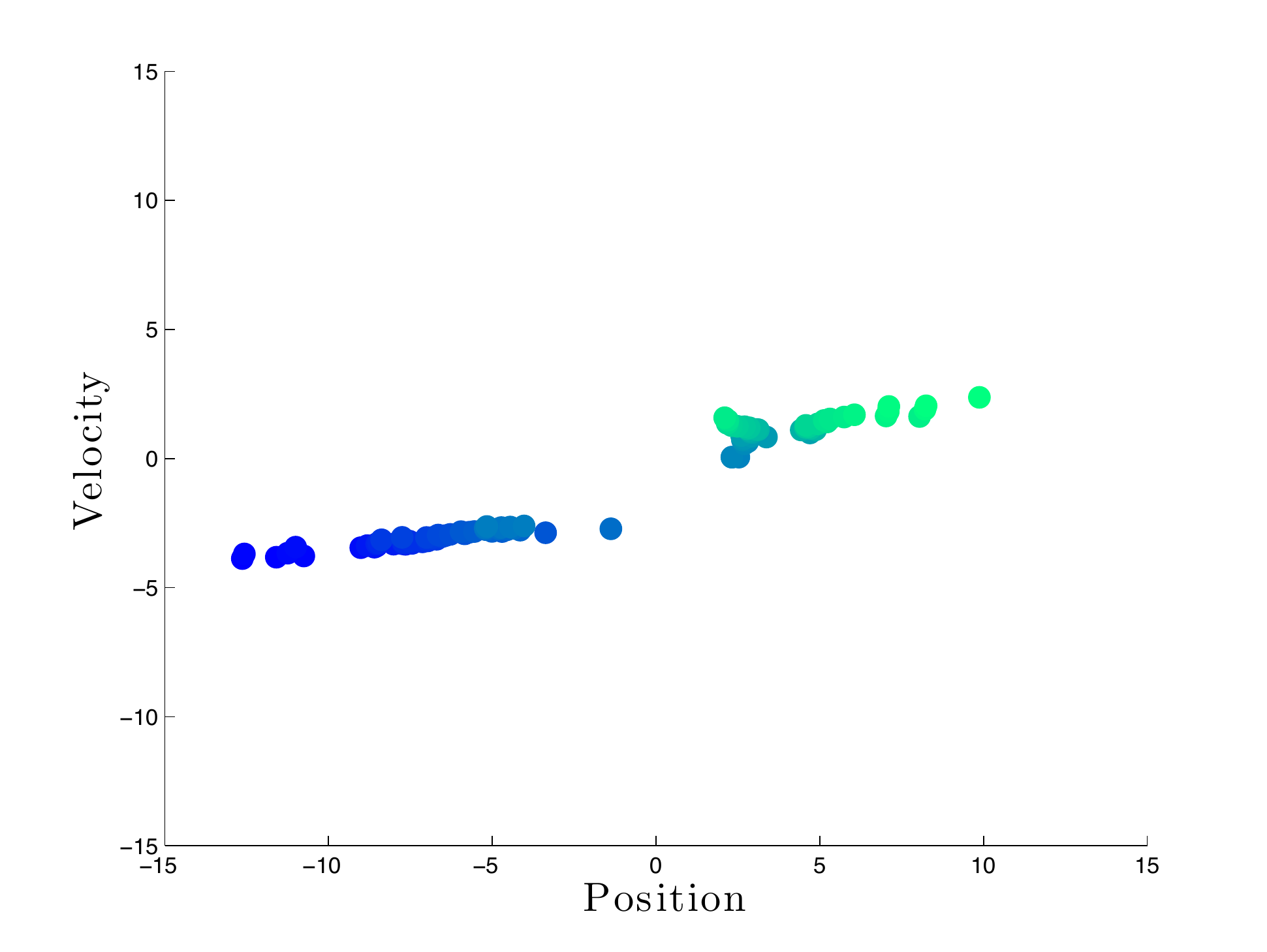}}
\hskip -0.7cm
\subfigure[t=5]{\includegraphics[scale=0.32]{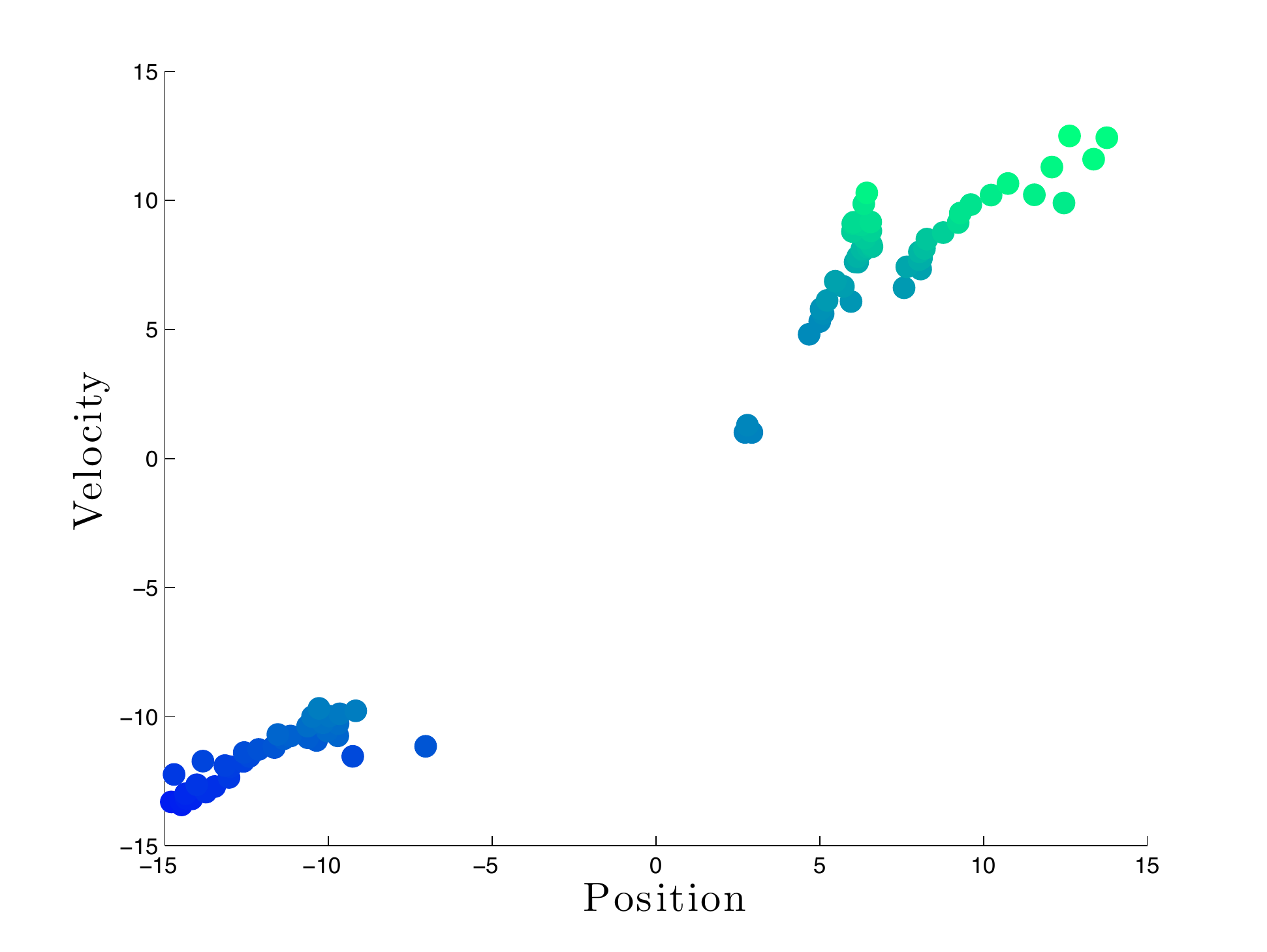}}\\
\caption{Numerical solution of \eqref{eq:CScontrolled}, with $\gamma=0.05<1/2$, $\zeta=0.01$, through a $10$th order gPC Hermite decomposition \eqref{eq:decomposition} with $\kappa=\infty$ with time step $\Delta t=10^{-2}$. The random input is normally distributed $\theta\sim\mathcal{N}(2,1)$. } \label{Fig:separation}
\end{center}\end{figure}

\begin{figure}[tbh]
\begin{center}
\subfigure[t=0]{\includegraphics[scale=0.3]{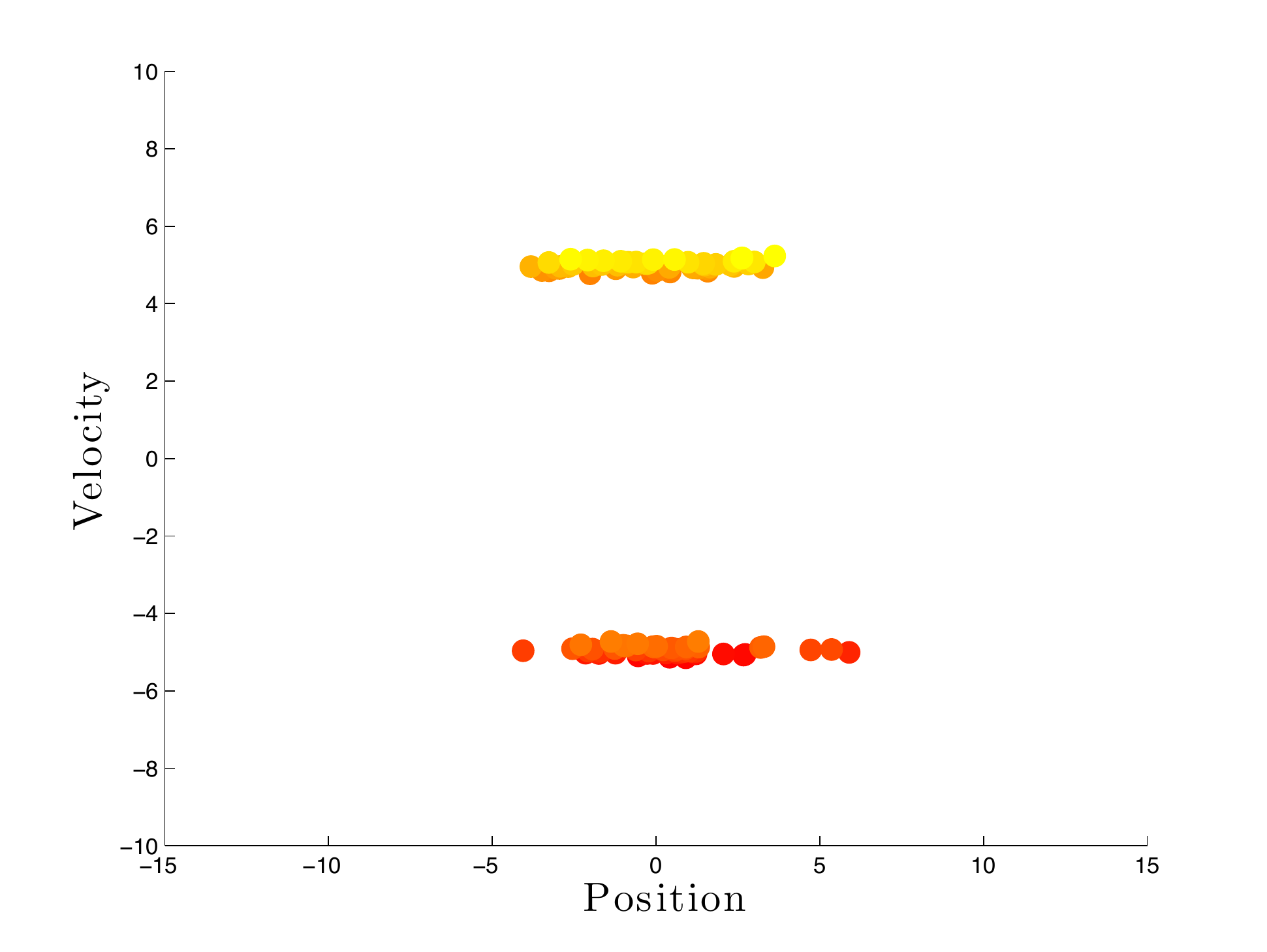}}
\subfigure[t=1]{\includegraphics[scale=0.3]{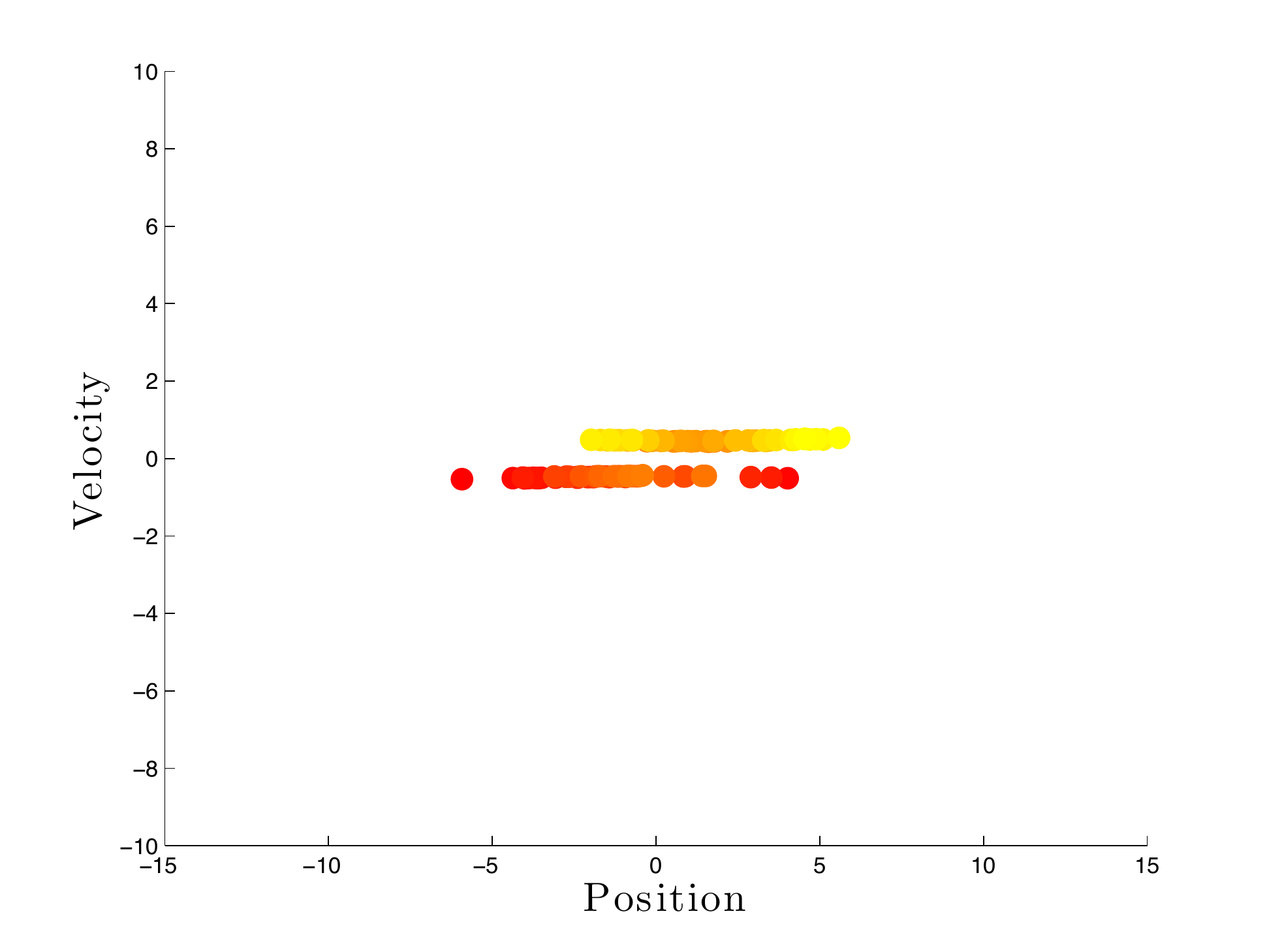}}\\
\subfigure[t=2]{\includegraphics[scale=0.3]{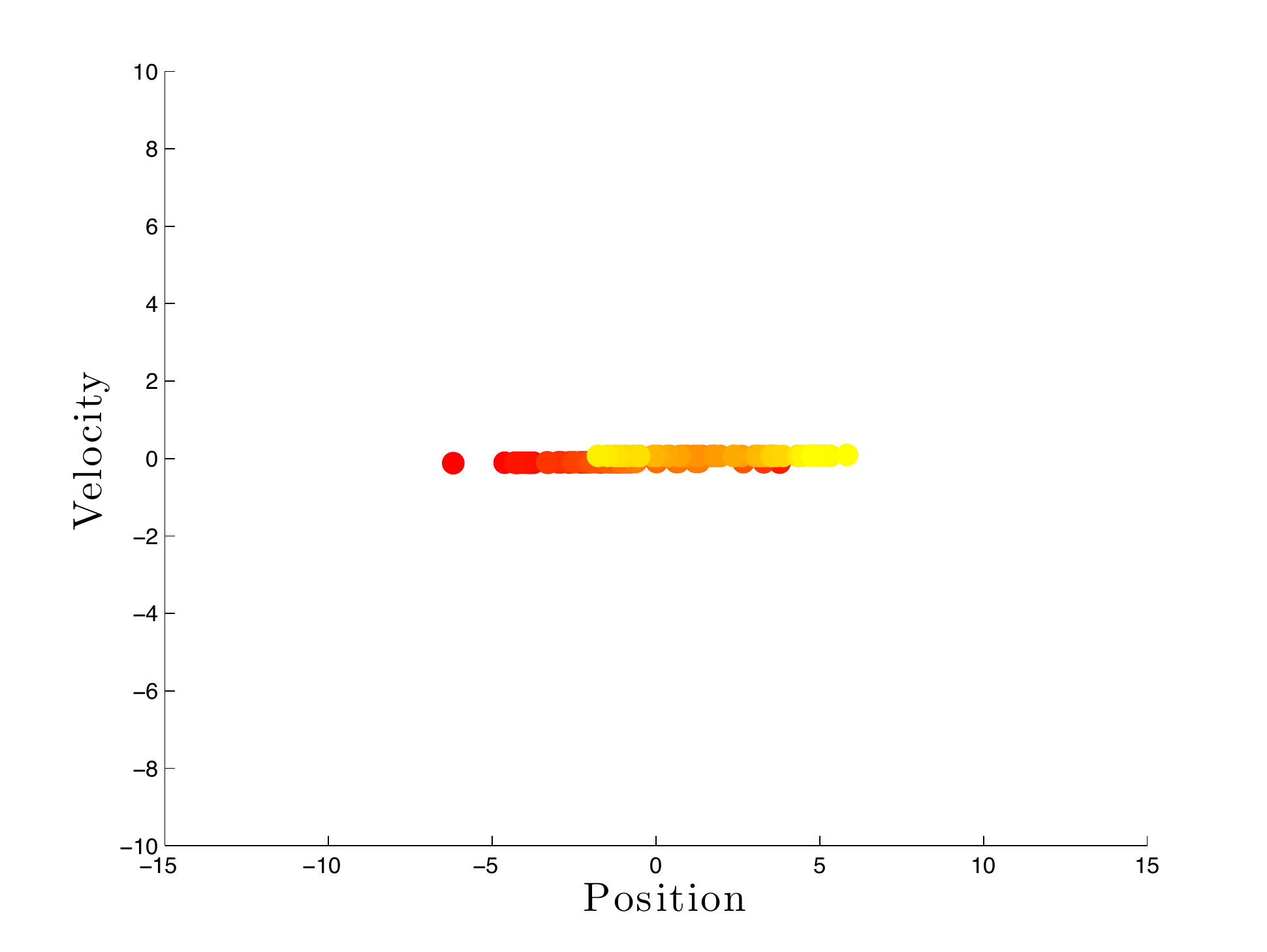}}
\subfigure[t=3]{\includegraphics[scale=0.3]{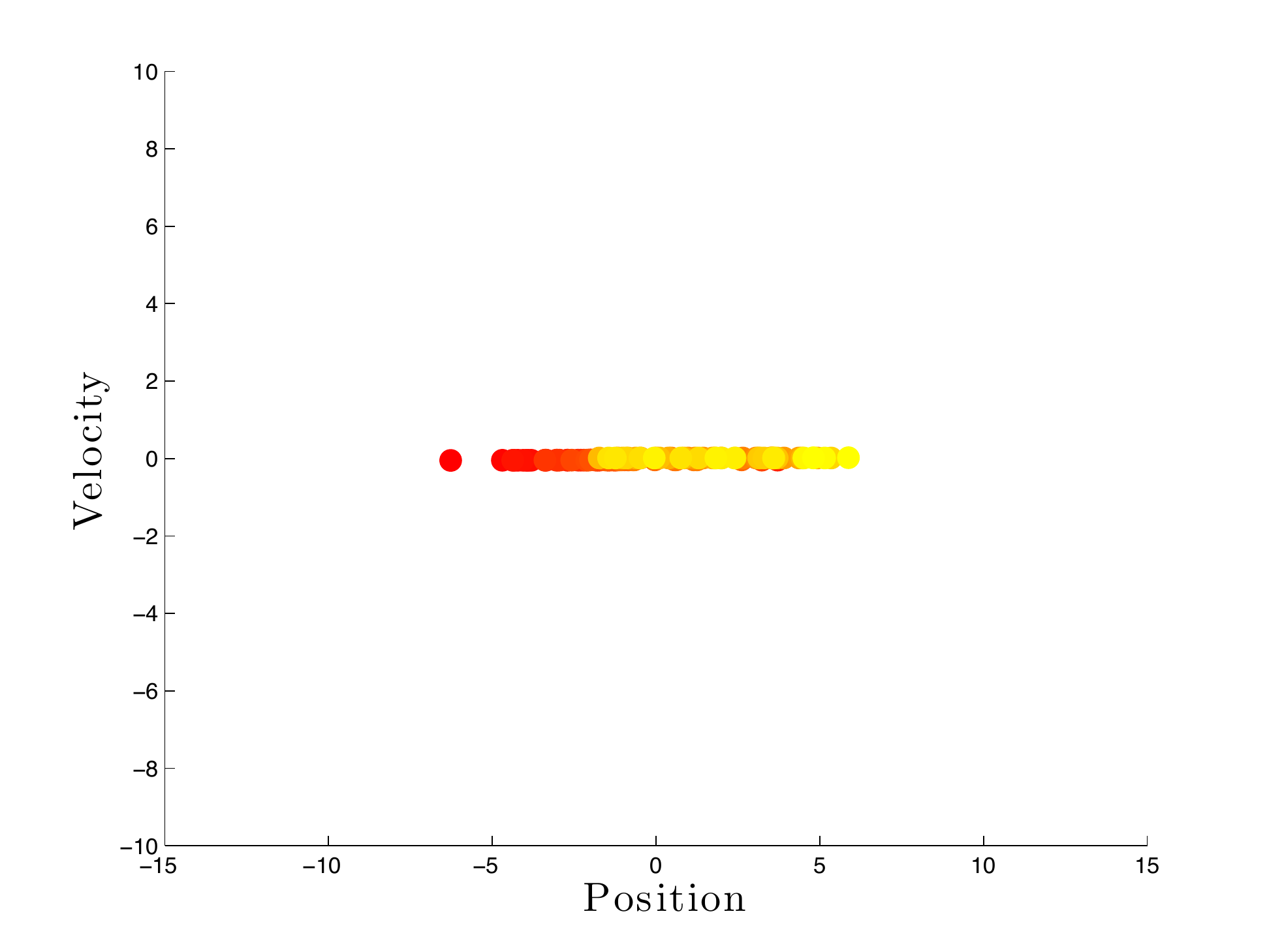}}
\subfigure[t=4]{\includegraphics[scale=0.3]{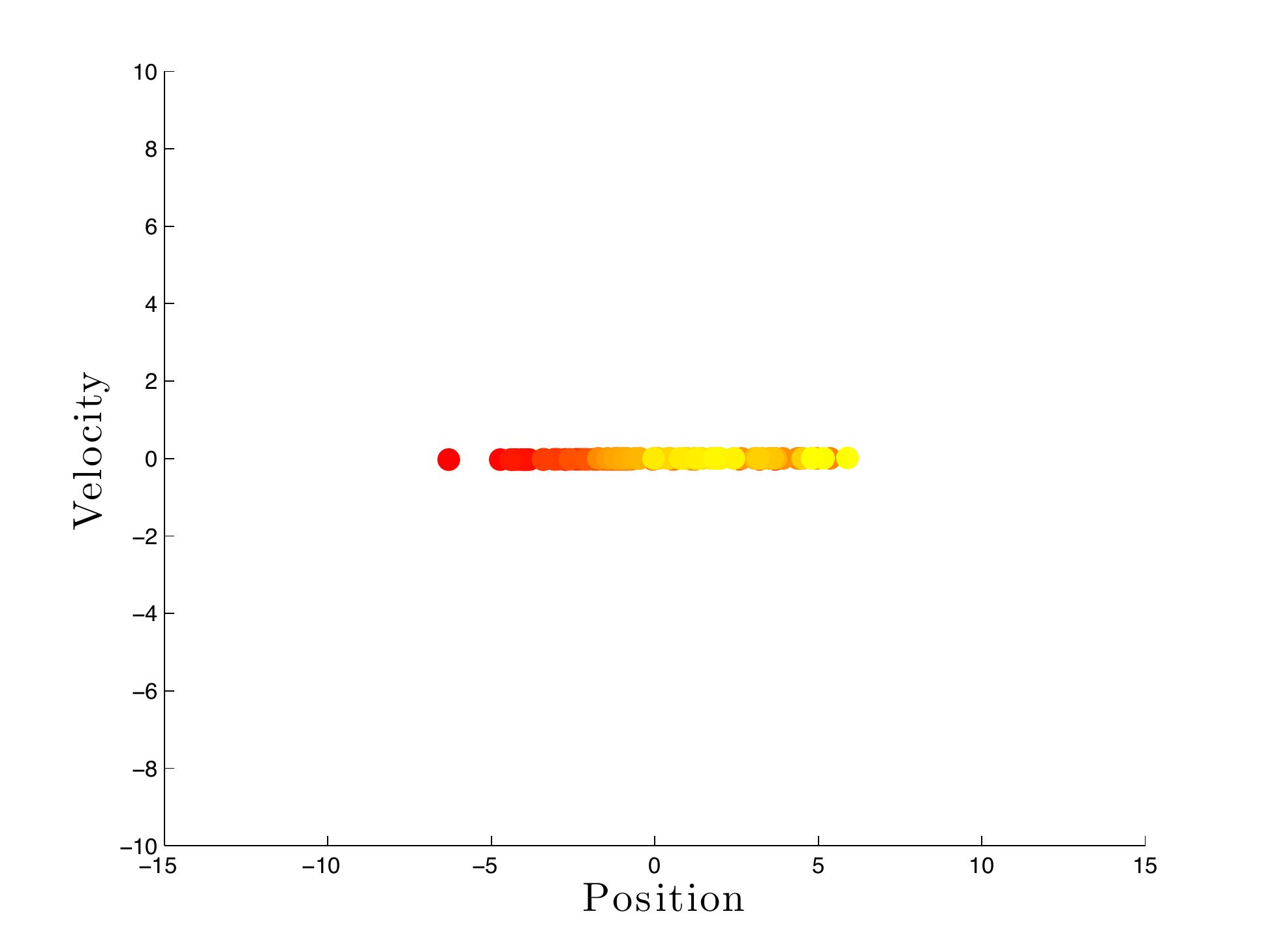}}
\subfigure[t=5]{\includegraphics[scale=0.3]{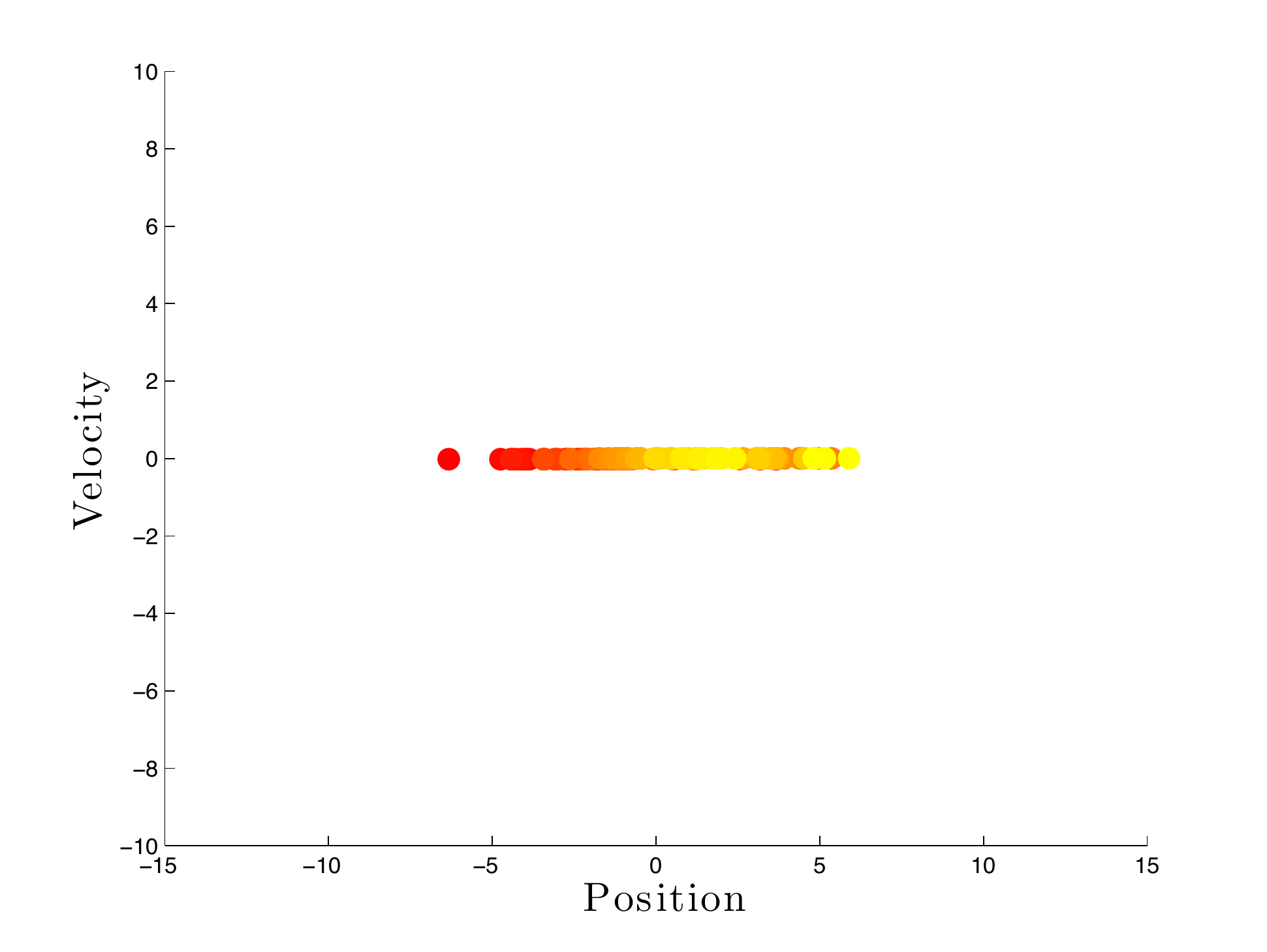}}
\caption{Numerical solution of \eqref{eq:CScontrolled}, with $\gamma=0.05<1/2$, $\zeta=0.01$, through a $6th$-order gPC Hermite decomposition for the selective control \eqref{eq:decomposition} with time step $\Delta t=10^{-2}$. Here we considered a normally distributed random input  $\theta\sim\mathcal{N}(2,1)$, the desired velocity is $v_d=0$ and the control parameter is $\kappa=1$.} \label{Fig:flocking}
\end{center}\end{figure}

\section{Conclusions}
\label{sec:Conc}
We proposed a general approach for the numerical approximation of flocking models with random inputs through gPC. The method is constructed in two steps. First the random Cucker-Smale system is solved by gPC. The presence of uncertainty in the interaction terms, which is a natural assumption in this kind of problems, leads to threshold effects in the asymptotic behavior of the system. Next a constrained gPC approximation is introduced and approximated though a selective model predictive control strategy. Relations under which the introduction of the gPC approximation and the model predictive control commute are also derived.
The numerical examples illustrates that the assumption of positivity of the mean value of the random input is not sufficient for the alignment of the system but that a suitable choice of the selective control is capable to stabilize the system towards the desired state. Extension of this technique to the case of a large number of interacting agents through mean-field and Boltzmann approximations are actually under study.

{

}
\end{document}